\documentclass[a4paper,12pt]{article}
\usepackage{amssymb}
\usepackage{lscape}
\usepackage{tabls}
\def\beq{\begin{equation}}
\def\eeq{\end{equation}}
\def\bea{\begin{eqnarray}}
\def\eea{\end{eqnarray}}
\def\nn{\nonumber}
\def\h#1{\hat{#1}}

\def\pair#1#2{\left\langle #1,#2 \right\rangle}
\def\O{OSp_q(1/2)}
\def\U{{\cal U}}
\def\A{{\cal A}}
\def\e#1#2{e^{#1}_{#2}}
\def\CGC#1#2{C^{#1}_{#2}\,}
\def\D#1#2{D^{\ #1}_{#2}}
\def\T#1#2{T^{\ #1}_{#2}}
\def\E#1#2{E^{#1}_{#2}(\Lambda)}
\def\oa{\mathfrak{a}}
\def\oab{\bar{\mathfrak{a}}}
\def\oc{\mathfrak{c}}
\def\bn#1#2{\left[\begin{array}{c} #1 \\ #2 \end{array}\right]}

\makeatletter
  \def\@cite#1#2{${\mbox{#1\if@tempswa , #2\fi}}$}
\makeatother
\def\ocite#1{$^{\mbox{\scriptsize{\cite{#1}}}}$}
\makeatletter
  \def\@biblabel#1{$^{\mbox{#1}}$}
\makeatother
\renewcommand{\theequation}{\arabic{section}.\arabic{equation}}
\renewcommand{\thesection}{\Roman{section}.}
\renewcommand{\thesubsection}{\Alph{subsection}.}
\begin{document}
%
%
%
%
%
%
\thispagestyle{empty}

\vspace*{3cm}
\begin{center}
{\Large\bf
Quantum Spheres for \boldmath{$ \O $}
}

\bigskip\bigskip
N. Aizawa

\bigskip
\textit{
Department of Applied Mathematics, \\
Osaka Women's University, \\
Sakai, Osaka 590-0035, Japan
}

\bigskip
and

\bigskip
R. Chakrabarti

\bigskip
\textit{
Department of Theoretical Physics, \\
University of Madras, \\
Guindy Campus, Chennai 600 025, India
}

\end{center}

\vfill
\begin{abstract}
  Using the corepresentation of the quantum supergroup $\O$ a general 
method for constructing noncommutative spaces covariant under its 
coaction is developed. In particular, a one-parameter family of covariant 
algebras, which may be interpreted as noncommutative superspheres, is 
constructed. It is observed that embedding of the supersphere 
in the $ \O $ algebra is possible. This realization admits 
infinitesimal characterization {\it \`a la} Koornwinder. A covariant 
oscillator realization of the supersphere is also presented. 
\end{abstract}
%
\newpage
\setcounter{page}{1}
\section{INTRODUCTION}

  Lie supergroups and superalgebras have been used as  basic tools 
in various fields of theoretical physics. 
Supersymmetry in quantum field theories and string field theory  
is the most well-known example of application of Lie superalgebras. 
Other examples are found in exactly solvable lattice models, interacting 
boson-fermion models in nuclear physics, extended t-J models in condensed 
matter physics, and so on. On the other hand, importance of noncommutative 
geometry in theoretical physics, especially in string theory and 
quantum gravity, has come into focus recently\ocite{Sz}. Therefore if 
these two notions are combined to form a noncommutative geometry of 
supersymmetric nature, we can expect that the combination  will play 
important roles in various fields in physics. An attempt to introduce a 
noncommutative superspace was made by Manin\ocite{Ma} in the context of 
quantum supergroup. Then differential calculus on the noncommutative 
superspace was developed by two different approaches\ocite{KU,So}. The 
present authors extended the notions of noncommutative differential 
geometry such as connection and curvature to the supersymmetric case and 
investigated\ocite{AC} the superspace for super-Jordanian deformed  
$OSp(1/2)$ group. 

 A supersphere having bosonic and fermionic coordinates is defined as 
an algebra whose defining relations are covariant under the coaction of 
the quantum supergroup $\O$. In this paper, we construct a one-parameter 
family of superspheres by  developing a {\it general procedure} based on 
the representation theory of $ \O $ and its dual $U_q[osp(1/2)]$. By 
this method, quantum superspaces and superspheres are described in a 
unified way. Furthermore, the method can  be used to find higher 
dimensional noncommutative superspaces covariant under $ \O.$ Our work 
is motivated by two reasons: (1) In order to investigate noncommutative 
geometry, it is important to have explicit examples of noncommutative 
superspaces. Manin's work is an analogue of flat space, while here we 
consider an analogue of curved 
space. (2) There exist some models of integrable quantum field theories 
with $ OSp(m/2n) $ symmetries where superspheres appear\ocite{SWK}. 
Consideration of quantum group extensions of such models will require 
quantum superspheres. We start with the simplest and the most important 
group $\O$ to construct quantum superspheres.   

  Let us briefly recall the studies of noncommutative sphere based on 
quantum groups, since supersymmetric counterparts of some of them are 
considered in this paper. Podle\'s introduced \ocite{Po1} an algebra 
which is covariant under the adjoint corepresentation of quantum group 
$ SU_q(2)$. The algebra is interpreted as a noncommutative 
version of two-sphere, and called $q$-sphere. The $q$-sphere has one 
more parameter in addition to its radius and the deformation parameter $q.$ 
Thus what Podle\'s constructed is a one-parameter family of noncommutative 
two-spheres. The parameter is specific to $q$-sphere and does not have 
commutative counterpart. Differential calculus on the $q$-sphere was 
initiated by Podle\'s\ocite{Po2}, then classification of differential 
structures on $q$-sphere was made in Refs. \cite{Po3}, \cite{AS}, 
\cite{HK}. An interesting relation of $q$-sphere to $q$-hypergeometric 
functions is discussed in Ref. \cite{NM}, where orthogonal bases on 
$q$-sphere are explicitly determined in terms of big $q$-Jacobi 
polynomials. The $q$-sphere can be realized by embedding it in $ SU_q(2).$ 
This embedding admits an elegant description \ocite{Koo} of $q$-sphere 
as an algebra which is invariant under left and right actions of a 
twisted primitive element of the quantum algebra $U_q[su(2)]$. Podle\'s 
$q$-sphere has been generalized in two different directions: higher 
dimensional, and Jordanian $SU(2). $ Higher dimensional $q$-spheres, 
more precisely, noncommutative analogue of $(2n+1)$-spheres were 
constructed in a similar way by replacing $ SU_q(2)$ with $ SU_q(n+1). $  
Furthermore, an invariant integral on quantum $ (2n+1)$-sphere is 
obtained in Ref. \cite{VS}. Another family of quantum two-spheres is 
obtained \ocite{CS} by using Jordanian 
deformation of $SL(2)$. One of its distinctions from Podle\'s $q$-sphere 
is that the Jordanian quantum sphere requires different twisted primitive 
elements for left and right invariances. 

  Throughout this paper, the quantum superalgebra $ U_q[osp(1/2)] $ and 
the quantum supergroup $ \O $ are denoted by $ \U $ and $ \A, $ respectively. 
We assume that $q$ is generic in this paper. The plan of this paper is as 
follows. In the first two preliminary sections we fix our notations 
and conventions, and list formulae used in 
subsequent sections. Section II is a summary of definitions and 
representation theory of $ \U. $ For computational purpose, 
we give all the defining relations of $ \A $ explicitly in Sec. \ref{QSG} 
A relation between representations of $ \U $ and corepresentations 
of $ \A $ is given in Sec. \ref{CorepA} A product law of two 
corepresentations, which is a quantum supergroup analogue of 
Wigner's product law in the quantum theory of angular momentum,  
is also derived in Sec. \ref{CorepA} A general prescription to find an 
algebra covariant under the coaction of $\A$ is given in 
Sec. \ref{Acovariant} As a simple application of the method, the most 
general form of quantum superspaces is derived. The method is applied to 
construct a one-parameter family of quantum superspheres in 
Sec. \ref{Asphere} Properties of the quantum supersphere are examined and 
the similarities to those for $q$-sphere are pointed out.

%
\section{\boldmath{$ U_q[osp(1/2)] $} AND ITS REPRESENTATIONS}
\label{Uqandrep}
\subsection{Definition and representations}
\label{DefReP}

The universal enveloping algebra $\U = U_q[osp(1/2)]$ is generated 
by the two even $ K^{\pm 1}$, and the two odd elements $ v_{\pm} $ 
satisfying the commutation properties \ocite{KR}
\bea
  & & K K^{-1} = K^{-1} K = 1, \qquad
      K v_{\pm} = q^{\pm 1/2} v_{\pm} K, \nn \\
  & & \{ v_+, v_- \} = -\frac{K^2-K^{-2}}{q^4-q^{-4}}. \label{Uqdef}
\eea
The Casimir element is given by
\beq
  C = \left( \frac{q^{1/2}K^2 - q^{-1/2}K^{-2}}{q^4-q^{-4}} \right)^2
    - \frac{qK^2 + q^{-1} K^{-2}}{(q+q^{-1})(q^2+q^{-2})}\,\,v_- v_+ 
    - (q^{1/2}+q^{-1/2})^2\,\,v_-^2 v_+^2.
  \label{ospCas}
\eeq
The coproduct $ (\Delta)$, the counit $ (\epsilon), $ and the antipode 
$ (S)$ maps read 
\bea
  & & \Delta(K^{\pm 1}) = K^{\pm 1} \otimes K^{\pm 1}, \qquad
      \Delta(v_{\pm}) = v_{\pm} \otimes K^{-1} + K \otimes v_{\pm},
      \label{coproUq} \\
  & & \epsilon(K^{\pm 1}) = 1, \qquad\qquad\qquad\quad \epsilon(v_{\pm}) = 0,
      \label{counitUq} \\
  & & S(K^{\pm 1}) = K^{\mp 1}, \qquad\qquad\quad\; S(v_{\pm}) = -q^{\mp 1/2} v_{\pm}.
      \label{SUq}
\eea

  The finite dimensional irreducible representations of the $ \U $ algebra 
are said to be of the grade-star \ocite{MM} type. Each irreducible 
representation is specified by a nonnegative integer $ \ell $ and the 
corresponding $(2\,\ell + 1)$ dimensional graded vector space $V^{(\ell)}$ 
that admits a nondegenerate Hermitian bilinear form denoted by $(\ ,\ )$. 
The subspaces of $ V^{(\ell)}$ having different parities are orthogonal with 
respect to the bilinear form. The graded adjoint operation $(*)$ is defined by
\beq
    (A^* f, g) = (-1)^{\h{A}\h{f}}(f,Ag), \qquad
    A \in \U;\;\; f, g \in V^{(\ell)},
    \label{GAO}
\eeq
where $\h{A} $ denotes the parity of $A$. The $*$-operation is assumed to be 
an algebra anti-isomorphism and coalgebra isomorphism: 
\beq
   (A_1 A_2)^* = (-1)^{\h{A_1}\h{A_2}} A_2^* A_1^*, \qquad
   (A_1 \otimes A_2)^* = A_1^* \otimes A_2^*.
   \label{starmorph}
\eeq 
The grade-star representation of $\U$ is characterized by
\beq
   K^* = K, \qquad v_{\pm}^* = \pm (-1)^{\varepsilon} v_{\mp},
   \qquad \varepsilon = 0, 1, 
   \label{GSRU}
\eeq
where $ \varepsilon $ refers to the class of the representation. 

  Let $ \{\; e^{\ell}_m(\lambda) \ | \ m = \ell, \ell-1, 
\cdots, -\ell \;\} $ be a basis of $ V^{(\ell)},$ where each basis 
vector has a definite parity. The index $ \lambda = 0, 1 $ 
specifies the parity of the highest weight vector $ e^{\ell}_{\ell}(\lambda). $ 
The parity of $ e^{\ell}_m(\lambda) $ equals $ \ell - m + \lambda,$ 
as it is obtained by the application of $ v_-^{\ell-m} $ on 
$ e^{\ell}_{\ell}(\lambda).$ For the superalgebras the norm of the 
representation basis need not be chosen positive definite.  
In this work, however, we assume the positive definiteness of 
the basis elements: 
\beq
 (e^{\ell}_m(\lambda), e^{\ell'}_{m'}(\lambda)) = \delta_{\ell \ell'} \delta_{mm'}.
 \label{Normalization}
\eeq
It turns out that this convention relates the parity $ \lambda $ 
and the class $ \varepsilon $ as follows:
\beq
  \lambda = \varepsilon + 1 \ ({\rm mod} \ 2).
  \label{lamandclass}
\eeq 
With these settings, the irreducible representations of $ \U $ are given by
\bea
  & & K \e{\ell}{m}(\lambda) = q^{m/2} \e{\ell}{m}(\lambda), \nn \\
  & & v_+ \e{\ell}{m}(\lambda) = \sqrt{[\ell-m][\ell+m+1]\varrho} 
      \;\e{\ell}{m+1}(\lambda),
      \label{RepU} \\
  & & v_- \e{\ell}{m}(\lambda) 
     = (-1)^{\ell-m-1} \sqrt{[\ell+m][\ell-m+1]\varrho}
     \; \e{\ell}{m-1}(\lambda), \nn
\eea
where $ [n] $ and $ \varrho $ are defined by
\beq
   [n] = \frac{q^{-n/2} - (-1)^n q^{n/2}}{ q^{-1/2}+ q^{1/2} } , \qquad
   \varrho = \frac{q^{-1/2} + q^{1/2}}{q^{-4}-q^4}.
   \label{Kulish}
\eeq
Our phase convention for $ v_{\pm} $ agrees with that of Ref. \cite{KR}, 
but it differs from that of Ref. \cite{MM}. For later convenience, the 
representation matrices for $ \ell = 1, 2 $ cases are given explicitly. 
The generators in the $ \ell=1 $ representation read
\bea
  & & K = {\rm diag} (q^{1/2},1,q^{-1/2}), \nn \\
  & & v_+ = \sqrt{[2]\varrho} 
          \left(
            \begin{array}{rrr}
               0 & 1 & 0 \\ 0 & 0 & 1 \\ 0 & 0 & 0
            \end{array}
          \right),
      \qquad
      v_- = \sqrt{[2]\varrho} 
          \left(
            \begin{array}{rrr}
               0 & 0 & 0 \\ -1 & 0 & 0 \\ 0 & 1 & 0
            \end{array}
          \right),
      \label{ell=1}  
\eea
and for the $ \ell = 2 $ case they are given by
\bea
  & & K = {\rm diag} (q,q^{1/2},1,q^{-1/2},q^{-1}), \nn \\
  & & v_+ = 
          \left(
            \begin{array}{rrrrr}
               0 & \sqrt{[4]\varrho} & 0 & 0 & 0 \\ 
               0 & 0 & \sqrt{[3]!\varrho} & 0 & 0 \\ 
               0 & 0 & 0 & \sqrt{[3]!\varrho} & 0 \\
               0 & 0 & 0 & 0 & \sqrt{[4]\varrho} \\
               0 & 0 & 0 & 0 & 0
            \end{array}
          \right),
      \label{ell=2} \\
  & & v_- = 
          \left(
            \begin{array}{rrrrr}
               0 & 0 & 0 & 0 & 0 \\
               -\sqrt{[4]\varrho} & 0 & 0 & 0 & 0 \\
               0 & \sqrt{[3]!\varrho} & 0 & 0 & 0 \\ 
               0 & 0 & -\sqrt{[3]!\varrho} & 0 & 0 \\
               0 & 0 & 0 & \sqrt{[4]\varrho} & 0
            \end{array}
          \right).
      \nn 
\eea
The eigenvalue of the Casimir element depends only on the highest weight:
\beq
   C \e{\ell}{m}(\lambda) = \left(\frac{q^{\ell+1/2} - q^{-\ell-1/2}}
   {q^4-q^{-4}} \right)^2\,\e{\ell}{m}(\lambda).
   \label{Cvalue}
\eeq

  The quantity $ [n]$, known\ocite{KR,MM} as Kulish symbol, plays a role 
similar to the $q$-number. For a positive integer $n$ its factorial  
is defined as $ [n]! \equiv [n] [n-1] \cdots [1].$ 

\subsection{Tensor product representations}
\label{TenRep}

  Tensor product of the irreducible representations of $\U$ has 
been discussed in Ref. \cite{KR}, \cite{MM}. 
The decomposition of tensor product in the irreducible representations 
is identical to the classical case
\[
  V^{(\ell_1)} \otimes V^{(\ell_2)} = 
  V^{(\ell_1+\ell_2)} \oplus V^{(\ell_1+\ell_2 -1)} \oplus \cdots 
  \oplus V^{(|\ell_1-\ell_2|)}.
\]
The explicit formulae of the Clebsch-Gordan coefficients (CGC) are 
obtained in Ref. \cite{MM}. We list below the relations which will be 
used in the later sections. 

  In spite of our assumption (\ref{Normalization}) of the positivity of 
the basis states, the norm of the tensor product of the bases is not 
always positive definite. For instance, the following norm can be negative 
for some combinations of $ \ell_a, m_a \ ( a = 1, 2) $ and $ \lambda$:
\bea
   & & (\e{\ell_1}{m_1}(\lambda) \otimes \e{\ell_2}{m_2}(\lambda),
   \e{\ell_1}{m_1}(\lambda) \otimes \e{\ell_2}{m_2}(\lambda)) 
   \nn \\
   & & \qquad \qquad
   = (-1)^{(\ell_1-m_1+\lambda)(\ell_2-m_2+\lambda)} (\e{\ell_1}{m_1}(\lambda),\e{\ell_1}{m_1}(\lambda))
   (\e{\ell_2}{m_2}(\lambda),\e{\ell_2}{m_2}(\lambda)).
   \label{normTen}
\eea
The irreducible basis of the tensor product representations is obtained 
by using the CGC:
\beq
    \e{\ell}{m}(\ell_1,\ell_2,\Lambda) = \sum_{m_1,m_2} 
    \CGC{\ell_1\ \ell_2\ \;\ell}{m_1\,m_2\,m}\,\,
    \e{\ell_1}{m_1}(\lambda) \otimes \e{\ell_2}{m_2}(\lambda),
    \label{CGCdef}
\eeq
where $m= m_1+m_2, $ and 
$ \Lambda = \ell_1+\ell_2+\ell \ ({\rm mod}\ 2) $ is the parity of the highest 
weight vector $ \e{\ell}{\ell}(\ell_1,\ell_2,\Lambda). $  
Since our phase convention for representations of $ v_{\pm} $ differs from 
that of Ref. \cite{MM}, we can not use the expression of CGC given 
therein. The CGC in our convention is explicitly given by
\bea
  & & \CGC{\ell_1\ \ell_2\ \;\ell}{m_1\,m_2\,m} 
  \nn \\
  & & \qquad = (-1)^{(\ell_1-\ell+m_2)(\ell-m+\lambda) + (\ell_1-\ell+m_2)(\ell_1-\ell+m_2+1)/2}
  \nn \\
  & & \qquad \times q^{m_2(m+1)/2 + (\ell_1-\ell_2)(\ell_1+\ell_2+1)/4-\ell(\ell+1)/4}
  \nn \\
  & & \qquad \times \left(
    [2\ell+1] \frac{[\ell_1+\ell_2-\ell]!\,[\ell+m]!\,[\ell-m]!\,[\ell_1-m_1]!\,[\ell_2-m_2]!}
                  {[\ell_1+\ell_2+\ell+1]!\,[\ell_1-\ell_2+\ell]!\,[-\ell_1+\ell_2+\ell]!\,
                  [\ell_1+m_1]!\,[\ell_2+m_2]!}
   \right)^{1/2}
  \nn\\
  & & \qquad \times
    \sum_{k}(-1)^{k(k-1)/2+ k(\ell_1+\ell_2-m)} q^{k(\ell+m+1)/2}
  \nn \\
  & & \qquad \times
      \frac{[\ell_1+\ell-m_2-k]!\,[\ell_2+m_2+k]!}
      {[k]!\,[\ell-m-k]!\,[\ell_1-\ell+m_2+k]!\,[\ell_2-m_2-k]!},
      \label{CGCformula}
\eea
where the index $k$ runs over all non-negative integers maintaining the 
arguments of $[x]$ non-negative. The derivation of (\ref{CGCformula}) is 
described in  \ref{CGCderivation}. All the CGC are, we note, of parity 
zero. The basis (\ref{CGCdef}) is pseudo orthogonal:
\beq
  (\e{\ell'}{m'}(\ell_1,\ell_2,\Lambda), \e{\ell}{m}(\ell_1,\ell_2,\Lambda))
  = 
  (-1)^{(\ell-m+\lambda)(\ell_1+\ell_2+\ell+\lambda)} 
  \delta_{\ell'\ell}\delta_{m'm}.
  \label{Normprod}
\eeq
The CGC satisfies two pseudo orthogonality relations
\beq
  \sum_{m_1,m_2} (-1)^{(\ell_1-m_1+\lambda)(\ell_2-m_2+\lambda)} 
  \CGC{\ell_1\ \ell_2\ \;\ell'}{m_1\,m_2\,m'} \CGC{\ell_1\ \ell_2\ \;\ell}{m_1\,m_2\,m} 
  =
  (-1)^{(\ell-m+\lambda)(\ell_1+\ell_2+\ell+\lambda)} \delta_{\ell\ell'}\delta_{mm'},
  \label{CGCortho1}
\eeq
\beq
  \sum_{\ell,m}(-1)^{(\ell-m+\lambda)(\ell_1+\ell_2+\ell+\lambda)}
  \CGC{\ell_1\ \ell_2\ \;\ell}{m_1'\,m_2'\,m} \CGC{\ell_1\ \ell_2\ \; \ell}{m_1\,m_2\,m}
  =
  (-1)^{(\ell_1-m_1+\lambda)(\ell_2-m_2+\lambda)} \delta_{m_1'm_1} \delta_{m_2'm_2}.
  \label{CGCortho2}
\eeq
Using (\ref{CGCortho2}), the construction (\ref{CGCdef}) is readily 
inverted:
\beq
   \e{\ell_1}{m_1}(\lambda) \otimes \e{\ell_2}{m_2}(\lambda) 
   = 
   (-1)^{(\ell_1-m_1)(\ell_2-m_2)} \sum_{\ell,m} (-1)^{(\ell-m)(\ell_1+\ell_2+\ell)}
   \CGC{\ell_1\ \ell_2\ \;\ell}{m_1\,m_2\,m} 
   \,\,\e{\ell}{m}(\ell_1,\ell_2,\Lambda).
   \label{CGCreverse}
\eeq

%
\setcounter{equation}{0}
\section{QUANTUM SUPERGROUP \boldmath{$\O $}}
\label{QSG}

  The quantum supergroup $ \A = \O $ is defined as a Hopf dual to the 
universal enveloping algebra $\U$ \ocite{KR}. In this section, all  
defining relations of $\A$ will be given explicitly. The universal 
${\cal R}$-matrix of $ \U $ is given in Ref. \cite{KR}. For the defining 
$ \ell=1$ representation, it reads
\beq
  {\cal R}^{(\ell = 1)} 
  = \left(
   \begin{array}{ccc|ccc|ccc}
      q & \cdot & \cdot & \cdot & \cdot & \cdot & \cdot & \cdot & \cdot \\
      \cdot & 1 & \cdot & \omega & \cdot & \cdot & \cdot & \cdot \\
      \cdot & \cdot & q^{-1} & \cdot & \lambda & \cdot & \rho & \cdot & \cdot \\
      \hline
      \cdot & \cdot & \cdot & 1 & \cdot & \cdot & \cdot & \cdot & \cdot \\
      \cdot & \cdot & \cdot & \cdot & 1 & \cdot & \lambda & \cdot & \cdot \\
      \cdot & \cdot & \cdot & \cdot & \cdot & 1 & \cdot & \omega & \cdot \\ 
      \hline
      \cdot & \cdot & \cdot & \cdot & \cdot & \cdot & q^{-1} & \cdot & \cdot \\
      \cdot & \cdot & \cdot & \cdot & \cdot & \cdot & \cdot & 1 & \cdot \\
      \cdot & \cdot & \cdot & \cdot & \cdot & \cdot & \cdot & \cdot & q
   \end{array}
  \right),
   \label{R-KR}
\eeq
where
\beq
 \omega = q-q^{-1}, \qquad \lambda = - q^{-1/2} \omega, \qquad 
 \rho = (1+q^{-1}) \omega,
 \label{param}
\eeq
and the dot is used instead of zero for better readability. 
Let $ P $ be a permutation operator: 
$ P(v \otimes w) = (-1)^{\hat{v}\hat{w}} w \otimes v. $ 
Standard FRT\ocite{FRT} construction is obtained via the matrix $R$:
\beq
 R = P {\cal R}^{(\ell = 1)} P = 
 \left(
   \begin{array}{rrr|rrr|rrr}
     q & \cdot & \cdot & \cdot & \cdot & \cdot & \cdot & \cdot & \cdot \\
     \cdot & 1 & \cdot & \cdot & \cdot & \cdot & \cdot & \cdot & \cdot \\
     \cdot & \cdot & q^{-1} & \cdot & \cdot & \cdot & \cdot & \cdot & \cdot \\
     \hline
     \cdot & \omega & \cdot & 1 & \cdot & \cdot & \cdot & \cdot & \cdot \\
     \cdot & \cdot & -\lambda & \cdot & 1 & \cdot & \cdot & \cdot & \cdot \\
     \cdot & \cdot & \cdot & \cdot & \cdot & 1 & \cdot & \cdot & \cdot \\
     \hline
     \cdot & \cdot & \rho & \cdot & -\lambda & \cdot & q^{-1} & \cdot & \cdot \\
     \cdot & \cdot & \cdot & \cdot & \cdot & \omega & \cdot & 1 & \cdot \\
     \cdot & \cdot & \cdot & \cdot & \cdot & \cdot & \cdot & \cdot & q
   \end{array}
 \right).
  \label{Rmat}
\eeq
The quantum $T$-matrix, whose elements generate the algebra $ \A $ is
given by
\beq
  T = (t_{ij}) = \left(
   \begin{array}{rrr}
      a & \alpha & b \\
      \gamma & e & \beta \\
      c & \delta & d
   \end{array}
  \right),
  \label{Tmatrix}
\eeq
where the entries in latin (greek) characters are of even (odd) parity. 
The parity of the supermatrix $T$ is zero, {\it i.e.}, 
$ \h{t_{ij}} = \h{i}+ \h{j}.$ 
The RTT-relation describes the exchange properties on the entries of $T$. 
The $q$-orthosymplectic condition reads 
\beq
    T^{st} C T = D C, \qquad T C^{-1} T^{st} = D C^{-1},
    \label{OS}
\eeq
where
\beq
  C = \left(
   \begin{array}{rrr}
     0 & 0 & -q^{-1/2} \\
     0 & 1 & 0 \\
     q^{1/2} & 0 & 0
   \end{array}
   \right),
  \qquad
  C^{-1} = \left(
   \begin{array}{rrr}
     0 & 0 & q^{-1/2} \\
     0 & 1 & 0 \\
     -q^{1/2} & 0 & 0
   \end{array}
   \right),
   \label{Cmat}
\eeq
and the super-determinant $D$ is given by
\beq
  D = ad - q bc - q^{1/2} \alpha \delta. 
  \label{sDet}
\eeq
The $ T^{st}$ denotes the super-transpose of $T$. The super-transpose of 
an arbitrary matrix is given as 
$ A^{st}_{ij} = (-1)^{\h{i} (\h{j}+1)} A_{ji}. $ The RTT-relations 
require $D$ to be central. 

The coproduct and counit of $T$ are given as usual
\beq
  \Delta(T) = T \stackrel{\cdot}{\otimes} T, \qquad
  \epsilon(T) = {\rm diag}(1,1,1).
  \label{coalgT}\eeq
The group-like property of the super-determinant 
$\Delta(D) = D \otimes D$ is obtained by taking the coproduct of both
sides of the relation (\ref{OS}). 
This allows us to set the constraint 
\beq 
D = ad - qbc - q^{1/2} \alpha \delta = 1. 
\label{det_1}
\eeq
The antipode of $T$ satisfies $ S(T) T = T S(T) = 1,$
and it explicitly reads 
\beq
 S(T) = 
 C^{-1} T^{st} C
 = \left(
  \begin{array}{ccc}
    d               & q^{-1/2} \beta  & -q^{-1} b \\
    -q^{1/2} \delta & e               & q^{-1/2} \alpha \\
    -qc             & -q^{1/2} \gamma & a
  \end{array}
  \right).
  \label{antipode}
\eeq

  The RTT-relations reveal that not all the entries of $T$ are independent. 
We express the elements $ e, \beta $ and $ \gamma $ in terms of the rest. 
Using (2,2) component of the first relation in (\ref{OS})
\beq
  e^2 = 1 - 2q^{-1/2} \alpha \delta + \omega bc,
  \label{esq1}
\eeq
we solve for $e$ and its inverse $e^{-1}$:
\begin{eqnarray}
 && e = 1 - q^{-1/2} \alpha \delta + (q-1) bc.\nonumber\\
 && e^{-1} = (1+q^{-1/2} \alpha \delta - (q-1) bc) 
    (1+q^{-1}(q-1)^2 bc)^{-1}.
  \label{einv}
\end{eqnarray}
Inclusion of the element $e^{-1}$ allows us to solve for $\beta$ and 
$\gamma$: 
\beq
  \beta = q^{-3/2} b \delta - q^{1/2} d \alpha, \qquad
  \gamma = q^{-1/2} a \delta - q^{3/2} c \alpha, \qquad
  \alpha \delta = \gamma \beta. \label{begam}
\eeq

  In summary, the algebra $\A$ is generated by $ a, b, c, d, \alpha $ and $ \delta.$ 
The generators satisfy the relations
\beq
 \begin{array}{lclcl}
   \vspace{1mm}
   ab = q^2 ba & \quad & ac = q^2 ca, & \quad & 
   [a,d] = -\lambda \alpha \delta +\rho bc, 
   \\ \vspace{1mm}
   a \alpha = q \alpha a, & & a \delta = q \delta a + q \omega c \alpha, & & 
   [b, c] = 0,
   \\ \vspace{1mm}
   bd = q^2 db, & & b \alpha = q^{-1} \alpha b, & & 
   b \delta = q \delta b,
   \\ \vspace{1mm}
   cd = q^2 dc, & & c \alpha = q^{-1} \alpha c, & & 
   c \delta = q \delta c,
   \\ \vspace{1mm}
   d \alpha = q^{-1} \alpha d - \omega \delta b, & & d \delta = q^{-1} \delta d, & & 
   \alpha \delta = -q \delta \alpha - q\lambda bc,
   \\ \vspace{1mm}
   \alpha^2 = -q^{-1}[2] ab, & & 
   \delta^2 = -q^{-1}[2]  cd. & & 
 \end{array}
 \label{OSPcom1}
\eeq
For later convenience, the commutation relations involving $ e, \beta $ 
and $ \gamma $ are listed below:
\beq
  \begin{array}{lclcl}
    \vspace{1mm}
    [a, e] = \omega \gamma \alpha, & \quad & a \beta = q \beta a + q\omega b \gamma, 
    & \quad & 
    a \gamma = q\gamma a,
    \\ \vspace{1mm}
    [b,e] = 0, & & b \beta = q \beta b, & & b\gamma = q^{-1} \gamma b, 
    \\ \vspace{1mm}
    [c,e] = 0, & & c \beta = q \beta c, & & c \gamma = q^{-1} \gamma c, 
    \\ \vspace{1mm}
    [d,e] = \omega \delta \beta, & & d \beta = q^{-1} \beta d, & & 
    d \gamma = q^{-1} \gamma d - \omega \beta c,
    \\ \vspace{1mm}
    [e, \alpha] = - \lambda \gamma b, & & [e, \beta] = \lambda b \delta, & & 
    [e, \gamma] = \lambda \alpha c,
    \\ \vspace{1mm}
    [e,\delta] = -\lambda c \beta, & & \{ \alpha, \beta \} = -\omega e b, & & 
    \{ \alpha, \gamma \} = 0,
    \\ \vspace{1mm}
    \beta \gamma = -q^{-1} \gamma \beta - \lambda bc, & & 
    \{ \beta, \delta \} = 0, & & \{ \gamma, \delta \} = \omega c e,
    \\ \vspace{1mm}
    \beta^2 = -q^{-1} [2] bd, & & 
    \gamma^2 = -q^{-1} [2] ac. & & 
  \end{array}
  \label{OSPcom2}
\eeq
Additional relations may also be proved:
\bea
  & & e \alpha = q^{1/2}(b \gamma - \beta a), \qquad\quad
      e \beta = q^{1/2}\delta b - q^{-1/2} \alpha d,
  \nn \\
  & & e \gamma = q^{1/2}(\delta a - c \alpha), \qquad\quad\;
      e \delta = q^{-1/2} \gamma d - q^{1/2} \beta c.
  \label{eodd}
\eea

  The duality between $ \U $ and $ \A $ is given by a nondegenerate 
pairing $ \pair{\ }{\ } $
\bea
  & & \pair{K}{T} = {\rm diag}(q^{1/2},1,q^{-1/2}), \nn \\
  & & \pair{v_{+}}{T} = \sqrt{[2]\varrho} 
      \left(
         \begin{array}{rrr}
            0 & 1 & 0 \\ 0 & 0 & -1 \\ 0 & 0 & 0
         \end{array}
      \right),
   \qquad
   \pair{v_-}{T} = \sqrt{[2]\varrho} 
      \left(
         \begin{array}{rrr}
            0 & 0 & 0 \\ 1 & 0 & 0 \\ 0 & 1 & 0
         \end{array}
      \right).
    \label{DualUqOq}
\eea
The pairing can be extended to tensor product algebras by setting 
\beq
  \pair{X_1 \otimes X_2}{a_1 \otimes a_2} = (-1)^{\h{X}_2\h{a}_1}\pair{X_1}{a_1} \pair{X_2}{a_2}.
  \label{TPpair}
\eeq

%
\setcounter{equation}{0}
\section{COREPRESENTATIONS OF $ \A $}
\label{CorepA}

  A vector space $V$ is called a right $\A$-comodule 
if there exists a linear mapping 
$ \varphi_R : V \rightarrow V \otimes \A $ satisfying
\beq
  (\varphi_R \otimes {\rm id})\circ \varphi_R = 
  ({\rm id}\otimes \Delta) \circ \varphi_R,
  \qquad
  ({\rm id}\otimes \epsilon) \circ \varphi_R = {\rm id}.
  \label{Rcomodule}
\eeq
Similarly, the left $\A$-comodule is defined as a vector 
space $V$ equipped with a linear mapping 
$ \varphi_L : V \rightarrow \A \otimes V $ such that
\beq
  ({\rm id}\otimes \varphi_L) \circ \varphi_L = 
  (\Delta \otimes {\rm id}) \circ \varphi_L,
  \qquad
  (\epsilon \otimes {\rm id}) \circ \varphi_L = {\rm id}.
  \label{Lcomodule}
\eeq
The mapping $ \varphi_R $ ($ \varphi_L $) is called 
a corepresentation, or, equivalently, a right (left) coaction of $ \A $ 
on $V$. 

  Employing the duality of the algebras $ \U $ and $ \A, $ we may 
follow the standard \ocite{KlS} construction of the action of $ \U $ 
on a $\A$-comodule $ V $. Namely, starting from the corepresentations of 
$ \A $, we may obtain representations of $ \U. $ Reversing the argument, 
we now obtain the hitherto unknown corepresentations of $ \A $ from the 
already known irreducible representations of $ \U. $  

  Let an arbitrary element $ X \in \U $ act on the vector space 
$ V^{(\ell)}$ described in \S \ref{Uqandrep}\ref{DefReP} 
The matrix representation of $ X$ on $ V^{(\ell)} $ is denoted by 
$ D^{\ell}(X;\lambda): $
\beq
   X \e{\ell}{m}(\lambda) = \sum_{m'} \e{\ell}{m'}(\lambda) \, 
   \D{\ell}{m'm}(X;\lambda).
   \label{Xrep}
\eeq
The parity of the element $ \T{\ell}{m'm}(\lambda) $ of $ \A $, defined 
via the duality relation 
\beq
   \D{\ell}{m'm}(X;\lambda) = 
   (-1)^{\h{X}(\ell-m'+\lambda)}\pair{X}{\T{\ell}{m'm}(\lambda)},
   \label{DtoT}
\eeq
may be assigned as 
\beq
 \widehat{\T{\ell}{m'm}}(\lambda) = m'+m \ ({\rm mod}\ 2). 
 \label{parityT} 
\eeq
To prove this, we compare the parity of the both sides of (\ref{Xrep}). 
The resultant equation 
$ \widehat{\D{\ell}{m'm}}(X;\lambda) = \h{X} + m'+m \ ({\rm mod}\ 2) $ 
establishes the relation (\ref{parityT}).  Though we know, via 
(\ref{RepU}), that all nonvanishing entries of the representation 
matrices are even, we should keep the parity of $ D^{\ell}(X;\lambda) $ 
in the present formulation. Next we obtain the coproduct and 
counit maps for $ \T{\ell}{m'm}(\lambda)$:
\beq
  \Delta(\T{\ell}{m'm}(\lambda)) = 
  \sum_{m''} \T{\ell}{m'm''}(\lambda) \otimes \T{\ell}{m''m}(\lambda),
  \qquad
  \epsilon(\T{\ell}{m'm}(\lambda)) = \delta_{m'm}.
  \label{DepsilonT}
\eeq
The matrix $ D^{\ell}(X;\lambda)$, being  a representation of 
$ X \in \U, $ obeys the rule
\beq
  \D{\ell}{m'm}(XY;\lambda) = \sum_{m''} \D{\ell}{m'm''}(X;\lambda) 
  \D{\ell}{m''m}(Y;\lambda) \qquad \forall X, Y \in \U. 
  \label{Dhom}
\eeq
Both sides of (\ref{Dhom}) may be rewritten by using (\ref{DtoT}):
\bea
 \rm{l.h.s.} &=& (-1)^{\widehat{XY}(\ell-m'+\lambda)}
 \pair{XY}{\T{\ell}{m'm}(\lambda)}\nn \\
 &=& (-1)^{\widehat{XY}(\ell-m'+\lambda)}
 \pair{X \otimes Y}{\Delta(\T{\ell}{m'm}(\lambda))},\nn \\
 \rm{r.h.s.} &=& \sum_{m''} (-1)^{\h{X}(\ell-m'+\lambda) + 
\h{Y}(\ell-m''+\lambda)} 
 \pair{X}{\T{\ell}{m'm''}(\lambda)} \pair{Y}{\T{\ell}{m''m}(\lambda)}
 \nn \\
 &=& \sum_{m''}(-1)^{(\h{X}+\h{Y})(\ell-m'+\lambda)} 
     \pair{X \otimes Y}{\T{\ell}{m'm''}(\lambda) \otimes \T{\ell}{m''m}(\lambda)}.
 \nn
\eea
The property (\ref{Dhom}) now assumes the form
\beq
  \pair{X \otimes Y}{\Delta(\T{\ell}{m'm}(\lambda)) 
     - \sum_{m''}\T{\ell}{m'm''}(\lambda) \otimes \T{\ell}{m''m}(\lambda)}
  =  0.
  \label{pairXYtoDT}
\eeq
Since the relation (\ref{pairXYtoDT}) is true for arbitrary $ X, Y $, and 
the pairing is nondegenerate, we obtain the coproduct map in 
(\ref{DepsilonT}). 
To obtain the counit map, we consider
\beq
  \D{\ell}{m'm}(1;\lambda) = \pair{1}{\T{\ell}{m'm}(\lambda)}
  \equiv \epsilon(\T{\ell}{m'm}(\lambda)), 
  \label{Dof1}
\eeq
and use the fact that the unit element of $\U$ is always represented by the 
identity matrix $(\delta_{m'm})$. This completes the proof of 
(\ref{DepsilonT}). 

   Defining the map $ \varphi_R : V^{(\ell)} \rightarrow V^{(\ell)} 
\otimes \A  $ by 
\beq
   \varphi_R(\e{\ell}{m}(\lambda)) = \sum_{m'} \e{\ell}{m'}(\lambda) \otimes 
   \T{\ell}{m'm}(\lambda),
   \label{RcoactionVl}
\eeq
it is easy to show, via (\ref{DepsilonT}), that $ V^{(\ell)} $ equipped 
with $ \varphi_R $ is a right $\A$-comodule. Thus the quantum supermatrix 
$ T^{\ell}(\lambda) $ provides the $(2\ell+1)$-dimensional corepresentation 
of $ \A $ on $ V^{(\ell)}.$  The relation between representations and 
corepresentations is summarized in the same form as non-super case\ocite{KlS}
\beq
   X \e{\ell}{m}(\lambda) = ({\rm id} \otimes X) \circ 
   \varphi_R(\e{\ell}{m}(\lambda)),
   \label{rep-corep}
\eeq
where the action of $ \U $ on $ \A $ is defined by the nondegenerate pairing. 
It is easy to find $T^{\ell}(\lambda)$ for $\ell = 0,\ 1$ from (\ref{DtoT}):
\beq
   \T{0}{00}(\lambda) = 1 \quad {\rm for} \ \lambda = 0, 1,
   \label{T000}
\eeq
\beq
  T^{1}(0) = T  = \left(
   \begin{array}{rrr}
      a & \alpha & b \\
      \gamma & e & \beta \\
      c & \delta & d
   \end{array}
  \right),
  \qquad
  T^{1}(1) = \left(
   \begin{array}{ccc}
      a & -\alpha & b \\
      -\gamma & e & -\beta \\
      c & -\delta & d
   \end{array}
  \right),
  \label{T1}
\eeq
where the indices of rows and columns of $ T^{1}(\lambda)$ run over $ 1, 0 $ 
and $-1.$ 

One of the important properties of the corepresentations $ T^{\ell}(\lambda)$ 
is that they satisfy the product law which gives a rule to combine two 
corepresentations to get the third one:
\bea
   & & \delta_{\ell'\ell} \T{\ell}{m'm}(\Lambda)
   = 
   \sum_{m_1,m_2 \atop m_1',m_2'} 
   (-1)^{p}
   \CGC{\ell_1\ \ell_2\ \;\ell'}{m_1'\,m_2'\,m'}
   \CGC{\ell_1\ \ell_2\ \;\ell}{m_1\,m_2\,m}
   \T{\ell_1}{m_1'm_1}(\lambda)\, \T{\ell_2}{m_2'm_2}(\lambda),
   \label{PLT} \\
   & & p = (m_1'+m_1)(\ell_2-m_2'+\lambda)+(\ell_1-m_1')(\ell_2-m_2')+
           (\ell'-m')(\ell_1+\ell_2+\ell').
   \nn
\eea
To prove the above product law, we first derive the fusion rule of 
representation matrices $ D^{\ell}(X;\lambda)$ of $ \U,$ and then 
apply the duality argument {\it \`a la} the proof of (\ref{DepsilonT}). 

Denoting the coproduct of $ X \in \U $ as 
$ {\displaystyle \Delta(X) = \sum_{a}X_a \otimes X^a, } $ 
we use the projection relation (\ref{CGCdef}) to obtain
\beq
  X\, \e{\ell}{m}(\ell_1,\ell_2,\Lambda) 
  = \sum_{m_1,m_2,a}  
  (-1)^{\h{X}^a (\ell_1-m_1+\lambda)} 
  \CGC{\ell_1\ \ell_2\ \;\ell}{m_1\,m_2\,m}\,\,
  X_a \e{\ell_1}{m_1}(\lambda) \otimes X^a \e{\ell_2}{m_2}(\lambda).
  \label{DXe}
\eeq
We employ (\ref{Xrep}) to compute both sides of (\ref{DXe}):
\bea
  \rm{r.h.s.} &=& \sum_{m_1,m_2 \atop m_1',m_2'} \sum_a  
  (-1)^{\h{X}^a (\ell_1-m_1+\lambda)}(-1)^{(\h{X}_a + m_1'+m_1)(\ell_2-m_2'+\lambda)}
    \nn \\
  & &  \qquad \qquad \times\,
  \CGC{\ell_1\ \ell_2\ \;\ell}{m_1\,m_2\,m}
  ( \e{\ell_1}{m_1'} \otimes \e{\ell_2}{m_2'} ) 
  \D{\ell_1}{m_1'm_1}(X_a;\lambda) \D{\ell_2}{m_2'm_2}(X^a;\lambda)
  \nn \\
  &=& \sum_{m_1,m_2 \atop m_1',m_2'} \sum_{a,\ell',m'}
  (-1)^{\h{X}^a (\ell_1-m_1+\lambda)+\h{X}_a (\ell_2-m_2'+\lambda)+p}
    \nn \\
  & & \quad \times\,
  \CGC{\ell_1\ \ell_2\ \;\ell'}{m_1'\,m_2'\,m'} \CGC{\ell_1\ \ell_2\ \;\ell}{m_1\,m_2\,m}
  \e{\ell'}{m'}(\ell_1,\ell_2,\Lambda^{\prime}) 
  \D{\ell_1}{m_1'm_1}(X_a;\lambda) \D{\ell_2}{m_2'm_2}(X^a;\lambda),
  \nn
\eea
where (\ref{CGCreverse}) is used in the last equality. 
The l.h.s. of (\ref{DXe}) follows directly from (\ref{Xrep}): 
\[
  {\rm l.h.s.} =  \sum_{m'} \e{\ell}{m'}(\ell_1,\ell_2,\Lambda)\, 
  \D{\ell}{m'm}(X;\Lambda), 
\]
yielding the product law of $ D^{\ell}(X,\lambda): $
\bea
  \delta_{\ell'\ell} \D{\ell}{m'm}(X,\Lambda)
  &=& 
  \sum_{m_1,m_2 \atop m_1',m_2'} \sum_{a}
  (-1)^{\h{X}^a (\ell_1-m_1+\lambda)+\h{X}_a(\ell_2-m_2'+\lambda)+p}
  \nn \\
  &\times& 
  \CGC{\ell_1\ \ell_2\ \;\ell'}{m_1'\,m_2'\,m'} \CGC{\ell_1\ \ell_2\ \;\ell}{m_1\,m_2\,m}
  \D{\ell_1}{m_1'm_1}(X_a;\lambda) \D{\ell_2}{m_2'm_2}(X^a;\lambda).
  \label{productD}
\eea
To derive the fusion rule for $ T^{\ell}(\lambda),$ we consider the 
dual pairing 
\bea
  & & \pair{X}{\T{\ell_1}{m_1'm_1}(\lambda)\, \T{\ell_2}{m_2'm_2}(\lambda)}
  \nn \\
  &=& \pair{\Delta(X)}{\T{\ell_1}{m_1'm_1}(\lambda) \otimes \T{\ell_2}{m_2'm_2}(\lambda)}
  \nn \\
  &=& \sum_{a} (-1)^{\h{X}^a(m'_1+m_1)}(-1)^{\h{X}_a(\ell_1-m_1'+\lambda)+\h{X}^a(\ell_2-m_2'+\lambda)}
  \D{\ell_1}{m_1'm_1}(X_a;\lambda) \D{\ell_2}{m_2'm_2}(X^a;\lambda),
  \nn
\eea
which, in turn, allows us to evaluate the following sum:
\bea
   & & \sum_{m_1,m_2 \atop m_1',m_2'} 
   (-1)^{p} \,
   \CGC{\ell_1\ \ell_2\ \;\ell'}{m_1'\,m_2'\,m'}
   \CGC{\ell_1\ \ell_2\ \;\ell}{m_1\,m_2\,m}
   \pair{X}{\T{\ell_1}{m_1'm_1}(\lambda)\, \T{\ell_2}{m_2'm_2}(\lambda)}
   \nn \\
   &=& (-1)^{\h{X}(\ell-m'+\Lambda)} 
    \sum_{m_1,m_2 \atop m_1',m_2'} \sum_{a}
  (-1)^{\h{X}^a (\ell_1-m_1+\lambda)+\h{X}_a(\ell_2-m_2'+\lambda)+p}
  \nn \\
  & &\times\,
  \CGC{\ell_1\ \ell_2\ \;\ell'}{m_1'\,m_2'\,m'} \CGC{\ell_1\ \ell_2\ \;\ell}{m_1\,m_2\,m}
  \D{\ell_1}{m_1'm_1}(X_a;\lambda) \D{\ell_2}{m_2'm_2}(X^a;\lambda)
  \nn \\
  &=& \delta_{\ell'\ell} (-1)^{\h{X}(\ell-m'+\Lambda)} \D{\ell}{m'm}(X;\Lambda)
  \nn \\
  &=&  \pair{X}{\delta_{\ell'\ell} \T{\ell}{m'm}(\Lambda)}.
  \nn
\eea
The second equality is due to (\ref{productD}). Invoking the arbitrariness 
of $X \in \U$, and the nondegeneracy of dual pairing $\pair{\ }{\ }$, we 
obtain the fusion rule (\ref{PLT}). 

%
\setcounter{equation}{0}
\section{$\A$-COVARIANT ALGEBRAS}
\label{Acovariant}

\subsection{General prescription}
\label{GP}

  In this section, we will give a general prescription to 
find $\A$-covariant algebras. By $\A$-covariant algebras, 
we mean algebras whose defining relations are covariant under 
the right coaction $ \varphi_R $ of $\A$ defined by (\ref{RcoactionVl}). 
Probably, the simplest way to 
find such an algebra is to introduce an algebraic structure on 
the representation space $ V^{(\ell)}.$ 
Assuming $\mu$ to be a multiplication map in $ V^{(\ell)},\ i.e.,
\ \mu(f \otimes g) = fg;\, f, g \in V^{(\ell)}, $ we specifically 
construct the following composite object:
\beq
  \E{L}{M} \equiv \mu( \e{L}{M}(\ell,\ell,\Lambda) )
  = \sum_{m_1,m_2} \CGC{\;\ell\ \;\,\ell\ \;\;L}{m_1\,m_2\,M} 
  \e{\ell}{m_1}(\lambda) \e{\ell}{m_2}(\lambda),
  \label{ELM}
\eeq
where $ \Lambda = L \ ({\rm mod}\ 2). $ 
The right coaction on $\E{L}{M}$ is shown to be
\beq
   \varphi_R(\E{L}{M}) = \sum_{M'} \E{L}{M'} \otimes \T{L}{M'M}(\Lambda).
   \label{phiELM}
\eeq
The proof may be done in a straightforward way by inverting the 
relation (\ref{ELM})
\beq
   \e{\ell}{m_1}(\lambda)  \e{\ell}{m_2}(\lambda) 
   = 
   (-1)^{(\ell-m_1)(\ell-m_2)} \sum_{L,M} (-1)^{(L-M)L}
   \CGC{\;\ell\ \;\ell\ \;\;L}{m_1\,m_2\,M} \E{L}{M},
   \label{ELMreverse}
\eeq
and subsequently using the product law (\ref{PLT})
\bea
 \varphi_R(\E{L}{M}) &=& 
    \sum_{m_1,m_2} \CGC{\;\ell\ \;\ell\ \;\;L}{m_1\,m_2\,M} 
    \varphi_R(\e{\ell}{m_1}(\lambda)) \varphi_R(\e{\ell}{m_2}(\lambda))
  \nn \\
  &=& \sum_{m_1,m_2 \atop m_1',m_2'} (-1)^{(m_1'+m_1)(\ell-m_2'+\lambda)}
      \CGC{\;\ell\ \;\ell\ \;\;L}{m_1\,m_2\,M} 
      \e{\ell}{m_1'}(\lambda) \e{\ell}{m_2'}(\lambda) \otimes 
      \T{\ell}{m_1'm_1}(\lambda) \T{\ell}{m_2'm_2}(\lambda)
  \nn \\
  &\stackrel{(\ref{ELMreverse})}{=}& 
       \sum_{m_1,m_2 \atop m_1',m_2'} \sum_{L',M'}(-1)^{p'}
      \CGC{\;\ell\ \;\ell\ \;\;L}{m_1\,m_2\,M} \CGC{\;\ell\ \;\ell\ \;\;L'}{m_1'\,m_2'\,M'} 
      \E{L'}{M'} \otimes 
      \T{\ell}{m_1'm_1}(\lambda) \T{\ell}{m_2'm_2}(\lambda)
  \nn \\
  &\stackrel{(\ref{PLT})}{=}&
      \sum_{M'} \E{L}{M'} \otimes \T{L}{M'M}(\Lambda).
  \nn
\eea
where $ p' = (m_1'+m_1)(\ell-m_2'+\lambda) + (\ell-m_1')(\ell-m_2') + (L'-M')L'. $ 

   Employing (\ref{phiELM}) we now extract a set of relations which are 
covariant under $\varphi_R.$ The $ L = 0 $ relation 
$ \varphi_R(E^{0}_{0}(0)) = E^{0}_{0}(0) $ signifies that  
$ E^{0}_{0}(0) $ is a scalar under the right coaction. It may be 
equated to a constant parameter $r$: 
\beq
  E^{0}_{0}(0) 
  = \sum_{m_1,m_2} \CGC{\;\ell\ \;\,\ell\ \;\;0}{m_1\,m_2\,0}\,\, 
  \e{\ell}{m_1}(\lambda) \e{\ell}{m_2}(\lambda)
  = r.
  \label{E00}
\eeq
If $ L = \ell $ and  $ \lambda = \ell \ ({\rm mod}\ 2), $ then 
$ E^{\ell}_m(\lambda) $ and $ \e{\ell}{m}(\lambda) $ have the same 
parity, and they transform identically under $\varphi_R$. Therefore 
$ E^{\ell}_m(\lambda) $ is, in general, proportional to 
$\e{\ell}{m}(\lambda)$. It may be noted that the following relations are 
covariant 
\beq
  E^{\ell}_m(\lambda)
  = \sum_{m_1,m_2} \CGC{\;\ell\ \;\,\ell\ \;\;\ell}{m_1\,m_2\,m}\,\, 
  \e{\ell}{m_1}(\lambda) \e{\ell}{m_2}(\lambda) 
  = \xi \e{\ell}{m}(\lambda),
  \label{E1m}
\eeq
where the proportionality constant $ \xi \rightarrow 0$ as $ q 
\rightarrow 1.$  For the case $ \lambda \neq \ell \ ({\rm mod}\ 2) $, 
$ E^{\ell}_m(\ell) $ has different parity from $ \e{\ell}{m}(\lambda),$ 
even though they transform identically. 
In this case, the constant $\xi $ in (\ref{E1m}) is regarded as a Grassmann 
number that also vanish at $ q = 1. $ For $L \neq 0, \ell$, the element  
$ E^{\ell}_m(\ell) $ can not be proportional to 
$\e{\ell}{m}(\lambda)$ as they transform differently. 
The relevant covariant relations are, therefore, of the form 
\beq
  \E{L}{M} = \sum_{m_1,m_2} \CGC{\;\ell\ \;\,\ell\ \;\;L}{m_1\,m_2\,M}\,\, 
  \e{\ell}{m_1}(\lambda) \e{\ell}{m_2}(\lambda)
  = 0.
  \label{ELMcov}
\eeq
As will be seen in the subsequent sections, the simultaneous use of all 
relations from (\ref{E00}) to (\ref{ELMcov}) gives an inconsistent result,  
since some of them do not have correct classical limits. In order to obtain 
a consistent covariant algebra, we have to make a choice regarding the
relations to be used for defining the algebra. Then the consistency has 
to be verified. As it is clear from the above discussion, the covariant 
algebras can have at most two more parameters ($r, \xi$) in addition to 
the deformation parameter $q.$ It is emphasised that the origin of the 
parameters is clearly explained in the framework of the representation theory. 

  We have formulated a method to construct $\A$-covariant algebras with 
respect to the right coaction. It is possible to repeat the same discussion 
for the left coaction. 

\subsection{Quantum superspace \boldmath{$(\ell = 1,\ \lambda = 0)$}}
\label{l10plane}

  Let us investigate the covariant algebra for $ \ell = 1$ case,  
where the relevant tensor product decomposition is given by $ 1 \otimes 
1 = 2 \oplus 1 \oplus 0. $ We assume that $ \lambda = 0 $, and denote 
the basis of $ V^{(1)} $ by $ z_m = \e{1}{m}(0).$ Thus $ z_{\pm 1} $ are 
parity even and $ z_0 $ is parity odd. The CGC for the decomposition is 
given in \ref{L1CGC}. For $ L = 0,$ we obtain from (\ref{E00})
\beq
   q^{1/2}z_{-1} z_1 + z_0^2 - q^{-1/2} z_1 z_{-1} = r.
   \label{rad10}
\eeq
For $ L = 1, $ we have $ \Lambda \neq \lambda $, and, therefore, the 
parameter $\xi$ is a Grassmann number: 
\bea
   & & -q^{1/2} z_0 z_1 + q^{-1/2} z_1 z_0 = \xi z_1,
   \nn \\
   & & z_{-1} z_1 + (q^{-1/2} + q^{1/2}) z_0^2 - z_1 z_{-1} = \xi z_0,
   \label{plane1} \\
   & & q^{1/2} z_{-1} z_0 - q^{-1/2} z_0 z_{-1} = \xi z_{-1}.
   \nn
\eea
For $ L = 2, $ we obtain, using (\ref{ELMcov}), unacceptable relations such as
\[
  z_1^2 = 0, \qquad q^{-1/2} z_0 z_1 + q^{1/2} z_1 z_0 = 0.
\]
Thus we take (\ref{rad10}) and (\ref{plane1}) as defining relations of 
our covariant algebra. 
We need to check the following conditions in order to verify 
whether or not the algebra is well-defined: 
\begin{enumerate}
\renewcommand{\labelenumi}{(\alph{enumi})}
 \item The constant $r$ commutes with all generators
 \item Product of three generators, say $ z_1 z_0 z_{-1}$, has two ways 
of reversing its ordering:
\[
  \begin{array}{ccccccc}
      &  & z_1 z_0 z_{-1} & \longrightarrow & z_1 z_{-1} z_0 & & \\
      & \nearrow & & & & \searrow &  \\
      z_0 z_1 z_{-1} & & & & & & z_{-1} z_1 z_0. \\
      & \searrow & & & & \nearrow & \\
      & & z_0 z_{-1} z_1 & \longrightarrow & z_{-1} z_0 z_1 & &
  \end{array}
\]
These two ways give the same result. 
\end{enumerate}
It is straightforward to verify that the condition (a) is satisfied. 
The condition (b), however, requires setting $ \xi = 0. $ 

  Therefore, we define our covariant algebra by combining relations 
(\ref{rad10}) and (\ref{plane1}), while maintaining $ \xi = 0$:
\bea
   & & z_1 z_0 = q z_0 z_1, \qquad z_0 z_{-1} = q z_{-1} z_0,
   \nn \\
   & & z_1 z_{-1} = q^2 z_{-1} z_1 - q(q^{-1/2}+q^{1/2}) r,
   \label{plane2} \\
   & & z_0^2 = -q^{-1}[2] z_1 z_{-1} - q^{-1}r.
   \nn 
\eea
This may be interpreted as the most general form of a quantum superspace. 
The simplest quantum superspace corresponds to the choice of $ r = 0. $

\subsection{Quantum superspace \boldmath{$(\ell = 1,\ \lambda = 1)$}}

   In this subsection, we study another quantum superspace where the 
parity is opposite to the previous example. Setting 
$ \theta_m = \e{1}{m}(1)$, we note that $ \theta_{\pm 1} $ have odd parity, 
while $ \theta_0 $ has even parity. Following  (\ref{E00}), we obtain
\beq
   q^{1/2} \theta_{-1} \theta_1 - \theta_0^2 - q^{-1/2} \theta_1 
   \theta_{-1} = r.
   \label{rad11}
\eeq 
The $ L = 1 $ relations may be obtained from (\ref{E1m}) except that 
in this case, as $\Lambda = \lambda$ holds, the parameter $\xi$ 
is not a Grassmann number. However, these relations such as
\[
   q^{1/2} \theta_0 \theta_1 + q^{-1/2} \theta_1 \theta_0 = \xi \theta_1.
\]
are unacceptable as they include, in the $q \rightarrow 1$ limit, 
anticommutators for the product of even and odd elements. 
On the other hand, $ L = 2 $ relations have proper classical limits:
\bea
  & & \theta_{\pm 1}^2 = 0, 
  \nn \\
  & & q^{-1/2} \theta_0 \theta_1 - q^{1/2} \theta_1 \theta_0 = 0,
  \label{plane3} \\
  & & q^{-1} \theta_{-1} \theta_1 - [2] \theta_0^2 + q \theta_1 \theta_{-1} = 0,
  \nn \\
  & & -q^{-1/2} \theta_{-1} \theta_0 + q^{1/2} \theta_0 \theta_{-1} = 0.
  \nn
\eea
An interesting observation for this case is that we have two kinds of 
quantum superspaces. Firstly we note that the relations (\ref{plane3}) 
are enough to define a covariant algebra, since it may be shown that the 
condition (b) of Sec. \ref{Acovariant}B is satisfied. Thus the relations 
(\ref{plane3}) define a quantum superspace. Alternately, combining 
(\ref{plane3}) with (\ref{rad11}) we obtain another set of 
covariant relations
\bea
   & & \theta_{\pm 1}^2 = 0, \qquad q \theta_1 \theta_0 = \theta_0 \theta_1,
   \qquad q \theta_0 \theta_{-1} = \theta_{-1} \theta_0,
   \label{plane4} \\
   & & \theta_1 \theta_{-1} + \theta_{-1} \theta_1 = -\frac{[2]}{[3]}r, 
   \qquad
   \theta_0^2 = -(q^{-1/2}+q^{1/2}) \theta_1 \theta_{-1} 
   - q^{1/2}\frac{[2]}{[3]}r,
   \nn
\eea
which have correct classical limit for an arbitrary value of $r$. 
With these relations, it may be checked that the conditions (a) 
and (b) of Sec \ref{Acovariant}B are satisfied. 
Thus we define the second quantum superspace by (\ref{plane4}). 

  The above quantum superspaces are covariant under the coaction
\beq
   \varphi_R(\theta_m) = \sum_{m'} \theta_{m'} \otimes \T{1}{m'm}(1),
   \label{Ractionplane}
\eeq
where $ T^1(1) $ is given by (\ref{T1}). Defining a new basis of the 
quantum superspace by
\beq
   \theta'_0 = \theta_0, \qquad \theta'_{\pm 1} = -\theta_{\pm 1},
   \label{newtheta}
\eeq
we observe that the $\theta'_m$ are covariant under the same 
corepresentation matrix as the case of $ \lambda = 0 $:
\beq
  \theta'_m = \sum_{m'} \theta'_{m'} \otimes \T{1}{m'm}(0). 
  \label{Ractionplane2}
\eeq

%
\setcounter{equation}{0}
\section{$\A$-COVARIANT SPHERE}
\label{Asphere}

\subsection{Construction}

  In this section, we investigate an algebra covariant under the coaction 
of $T^2(\lambda),$ $i.e.$ the adjoint corepresentation of $ \A. $ 
This may be interpreted as a supersymmetric extension of a noncommutative 
sphere. The corepresentation matrices $ T^2(\lambda) $ are found 
from (\ref{ell=2}) and (\ref{DtoT}), or, alternately,  by coupling the 
elements of two $ T^1(\lambda)$ matrices via the product law (\ref{PLT}). 
We restrict ourselves to the case of $ \lambda = 0: $ 
\beq
  T^2(0) = \left(
     \begin{array}{ccccc}
        a^2 & \kappa_1 a \alpha & \kappa_3 ab & \kappa_1 \alpha b & b^2
        \\
        \kappa_1 a \gamma & ae + q^{-1}\gamma \alpha & \kappa_2(a\beta + q^{-1} \gamma b)
        & -\alpha \beta + q^{-1} e b & \kappa_1 b \beta
        \\
        \kappa_3 ac & \kappa_2(a \delta + c \alpha) 
        & ad + q^{-1}[2] \alpha \delta + q^{-2} b c & \kappa_2(\alpha d + \delta b) 
        & \kappa_3 bd
        \\
        \kappa_1 \gamma c & \gamma \delta + q^{-1} c e 
        & \kappa_2(\gamma d + q^{-1} c \beta) & ed + q^{-1}\beta \delta
        & \kappa_1 \beta d
        \\
        c^2 & \kappa_1 c \delta & \kappa_3 c d & \kappa_1 \delta d & d^2
     \end{array}
  \right),
  \label{T2mat0}
\eeq
where
\beq
   \kappa_1 = \sqrt{\frac{[4]}{q[2]}} 
   \qquad
   \kappa_2 = \sqrt{q^{-1}[3]} 
   \qquad
   \kappa_3 = \kappa_1 \kappa_2.
   \label{kappas}
\eeq
  
  The basis of $ V^{(2)} $ is denoted by $ Y_m = \e{2}{m}(0)$, where   
$m = 0, \pm 1, \pm 2.$ Here $ Y_{0}, Y_{\pm 2} $ are even, 
and $ Y_{\pm 1} $ are odd. Following the prescription in 
Sec.\ref{Acovariant}A, we seek a covariant algebra under the right 
coaction of $ T^2(0). $ The CGC for $ \ell = 2 $ are found in Appendix D. 
The relation for $ L = 0 $ is obtained via (\ref{E00}):
\beq
  q^{-1} Y_2 Y_{-2} - q^{-1/2} Y_1Y_{-1} - Y_0^2 + q^{1/2} Y_{-1} Y_1 
  + q Y_{-2} Y_2 = r, \label{radius1}
\eeq
where $ r$ is a constant. Equation (\ref{radius1}) may be regarded as 
the radius relation of the quantum supersphere. Explicit constructions 
for  the $ L = 2 $ case are obtained from (\ref{E1m}):
\bea
  & & q^{-3/2} Y_2 Y_0 - \left( \frac{[3]!}{[4]} \right)^{1/2} 
  Y_1^2 - q^{3/2} Y_0 Y_2 = \xi Y_2,\nn \\
  & & q^{-1/2} \left( \frac{[3]!}{[4]} \right)^{1/2} Y_2 Y_{-1} + q^{-1/2} 
  \mu_{21} Y_1 Y_0 - q^{1/2} \mu_{22} Y_0 Y_1 - q^{1/2} 
  \left( \frac{[3]!}{[4]} \right)^{1/2} Y_{-1} Y_2 = \xi Y_1,\nn \\
  & & q^{1/2} Y_2 Y_{-2} - \mu_{22} Y_1 Y_{-1} + \mu_{23} Y_0^2 - \mu_{21} 
  Y_{-1} Y_1 - q^{-1/2} Y_{-2} Y_2 = \xi Y_0,
  \label{sphere1-1} \\
  & & q^{-1/2} \left( \frac{[3]!}{[4]} \right)^{1/2} Y_1 Y_{-2} + q^{-1/2} 
  \mu_{21} Y_0 Y_{-1} - q^{1/2} \mu_{22} Y_{-1} Y_0 - q^{1/2} 
  \left( \frac{[3]!}{[4]} \right)^{1/2} Y_{-2} Y_1 = \xi Y_{-1},\nn \\
  & & q^{-3/2} Y_0 Y_{-2} - \left( \frac{[3]!}{[4]} \right)^{1/2} Y_{-1}^2 
  - q^{3/2} Y_{-2} Y_0 = \xi Y_{-2},
  \nn
\eea
where $ \xi $ is a constant vanishing in the classical limit, and 
\beq
  \mu_{21} = \frac{[2]}{[4]} (q^{-2} + q^{3/2}[2]), \quad
  \mu_{22} = \frac{[2]}{[4]} (q^2-q^{-3/2}[2]), \quad
  \mu_{23} = \frac{[2]^2}{[4]} \left( \frac{[8]}{[4]} + 2\frac{[4]}{[2]} + 1 \right).
  \label{mu2}
\eeq
The construction (\ref{ELMcov}) may be applied to the remaining values 
of $L.$ The relations for the $ L = 3$ case read  
\bea
  & & q^{-1} Y_2 Y_1 - q Y_1 Y_2 = 0, 
  \nn \\
  & & Y_2 Y_0 - \mu_{31} Y_1^2 - Y_0 Y_2 = 0,
  \nn \\
  & & q Y_2 Y_{-1} - \mu_{32} Y_1 Y_0 + \mu_{33} Y_0 Y_1 - q^{-1} Y_{-1} Y_{2} = 0, 
  \label{sphere1-2} \\
  & & -q^2 Y_2 Y_{-2} + q^{1/2}([3]+q^2) Y_1 Y_{-1} + [3] \omega Y_0^2 
  + q^{-1/2} ([3]+q^{-2}) Y_{-1} Y_1 + q^{-2} Y_{-2} Y_2 = 0,
  \nn \\
  & & q Y_1 Y_{-2} - \mu_{32} Y_0 Y_{-1} + \mu_{33} Y_{-1} Y_0 - q^{-1} Y_{-2} Y_1 = 0,
  \nn \\
  & & Y_0 Y_{-2} - \mu_{31} Y_{-1}^2 - Y_{-2} Y_0 = 0,
  \nn \\
  & & q^{-1} Y_{-1} Y_{-2} - qY_{-2} Y_{-1} = 0,
  \nn
\eea
where
\beq
 \mu_{31} = -\frac{\sqrt{ [4]! }}{[2]^2} \omega,\quad
 \mu_{32} = \left( \frac{[4]}{[3]!} \right)^{1/2} (q^2+q^{-1/2}[2]), \quad
 \mu_{33} = \left( \frac{[4]}{[3]!} \right)^{1/2} (q^{-2}-q^{1/2}[2]).
 \label{mu3}
\eeq
We observe that the relations (\ref{sphere1-1}) together with 
(\ref{sphere1-2}) have the correct classical limits. We have obtained a 
set of twelve relations for five generators. To test whether they  
consistently define an algebra, we need to check for the conditions (a) 
and (b) mentioned in Sec.\ref{Acovariant}\ref{l10plane} 
It may be proved by direct computation that the said conditions 
are, however,  not satisfied. 

  In order to make the algebra well-defined, we incorporate the $ L=1 $ 
relations listed below:
\bea
  & & q^{-3/2} Y_2 Y_{-1} - q^{-1/2} \left( \frac{[3]!}{[4]} \right)^{1/2} Y_1 Y_0 
  - q^{1/2} \left( \frac{[3]!}{[4]} \right)^{1/2} Y_0 Y_1 + q^{3/2} Y_{-1} Y_2 = 0,
  \nn \\
  & & -q^{-1/2} Y_2 Y_{-2} + q^{-1} \mu_{11} Y_1 Y_{-1} + 
       \varpi \frac{[3]!}{[4]} Y_0^2 - q\mu_{12} Y_{-1} Y_1 - q^{1/2} Y_{-2} Y_2 = 0,
  \label{sphere1-4} \\
  & & -q^{-3/2} Y_1 Y_{-2} + q^{-1/2} \left( \frac{[3]!}{[4]} \right)^{1/2} Y_0 Y_{-1} 
  + q^{1/2} \left( \frac{[3]!}{[4]} \right)^{1/2} Y_{-1} Y_0 - q^{3/2} Y_{-2} Y_1 = 0,
  \nn
\eea
where
\beq
  \varpi = q^{1/2}+q^{-1/2}, \qquad \mu_{11} = q + q^{-1}\frac{[2]}{[4]}, \qquad
  \mu_{12} = q^{-1} + q \frac{[2]}{[4]}.
  \label{mu1} 
\eeq
The remaining $ L = 4$ relations can not be incorporated, 
because they contain unacceptable equations such as
$
  Y_{\pm 2}^2 = 0. 
$

   As all the relations in (\ref{sphere1-1}), (\ref{sphere1-2}), and 
(\ref{sphere1-4}) are covariant by construction, their  linear combinations 
are also covariant. Taking linear combination of the fifteen relations, twelve "commutation relations" of generators and 
three constraints are obtained. The commutation relations read
\bea
  & & Y_2 Y_1 = q^2 Y_1 Y_2, \qquad 
      Y_{-1} Y_{-2} = q^2 Y_{-2} Y_{-1}, \nn \\
  & & q^{-2} Y_2 Y_0 = q^2 Y_0 Y_2 + \varpi \frac{[4][3]}{[6]}\xi Y_2, \nn \\
  & & q^{-2} Y_0 Y_{-2} = q^2 Y_{-2} Y_0 + \varpi \frac{[4][3]}{[6]}\xi Y_{-2}, \nn \\
  & & q^{-3} Y_2 Y_{-1} = q^3 Y_{-1} Y_2 + \varpi \frac{[3]\sqrt{[4]!}}{[6]}\xi Y_1, \nn \\
  & & q^{-3} Y_1 Y_{-2} = q^3 Y_{-2} Y_1 + \varpi \frac{[3]\sqrt{[4]!}}{[6]}\xi Y_{-1}, 
      \label{sphere2-1} \\
  & & q^{-1} Y_1 Y_0 = q Y_0 Y_{1} + \varpi \frac{[3]!}{[6]}\xi Y_1, \nn \\
  & & q^{-1} Y_0 Y_{-1} = q Y_{-1} Y_0 + \varpi \frac{[3]!}{[6]}\xi Y_{-1}, \nn \\
  & & q^{-1} \mu_{11} Y_2 Y_{-2} = q \mu_{12} Y_{-2} Y_2 - F_1 Y_0^2 + F_2 Y_0, \nn \\
  & & q^{-1/2} Y_1 Y_{-1} = - q^{1/2} Y_{-1} Y_1 + \omega Y_0^2 + \frac{\varpi}{q+1+q^{-1}}\xi Y_0,
  \nn \\
  & & Y_1^2 = \frac{1}{\mu_{31}}(Y_2 Y_0 - Y_0 Y_2),
      \qquad
      Y_{-1}^2 = \frac{1}{\mu_{31}}(Y_0 Y_{-2} - Y_{-2} Y_0),
  \nn
\eea
where
\bea
  & & F_1 = \varpi (q^2+q+1+q^{-1}+q^{-2}) \frac{[2]^2}{[4]}, \nn \\
  & & F_2 = \varpi (q^2+2q+2q^{-1}+q^{-2}) \frac{[3]![2]}{[6][4]}\xi. \label{sphere2-2}
\eea
Classical commutation properties immediately follow from 
(\ref{sphere2-1}) in the limit $ q \rightarrow 1.$ 
Three said constraints are given by
\bea
  & & Y_2 Y_{-1} = \frac{q^2}{[3]} \left(\frac{[3]!}{[4]} \right)^{1/2} Y_1 Y_0
     + q^{1/2}\frac{[2]}{[6]} \left(\frac{[3]!}{[4]} \right)^{1/2} \xi Y_1,
     \nn \\
  & & Y_1 Y_{-2} = \frac{q^2}{[3]} \left(\frac{[3]!}{[4]} \right)^{1/2} Y_0 Y_{-1} 
     + q^{1/2}\frac{[2]}{[6]} \left(\frac{[3]!}{[4]} \right)^{1/2} \xi Y_{-1},
     \label{sphere2-3} \\
  & & Y_0^2 = q^{-1}\frac{[4]}{[2]} Y_2 Y_{-2} - q^{-1/2}(q+q^{-1}) \mu_{12} Y_1 Y_{-1}
      - q^{-3/2} \frac{[3]!}{[6]} \xi Y_0.
     \nn
\eea
The classical limit of these constraints are not required in the 
commutative case. However, {\it we need the constraints to satisfy the two 
important requirements} given in Sec.\ref{Acovariant}\ref{l10plane}  
Verification of these conditions is straightforward but requires dull 
lengthy computation. Since the two conditions are satisfied, we define 
a one-parameter family of $ \A$-covariant quantum superspheres  
by the radius relation (\ref{radius1}), the commutation relations 
(\ref{sphere2-1}) and the constraints (\ref{sphere2-3}). We denote this 
quantum supersphere by $ S_{q,\xi}^0. $ 
The superscript $0$ indicates the parity  $ \lambda = 0 $ and $ \xi $ is 
a free parameter 
which does not have classical counterpart. The origin of the parameter and 
the fact that 
the quantum superspheres can not have more parameters clearly follows from 
the formulation in Sec. \ref{Acovariant}\ref{GP}

\subsection{Properties of \boldmath{$ S_{q,\xi}^0$}}

  In this subsection, three properties of the quantum supersphere 
$ S_{q,\xi}^0 $ are investigated. First, we consider a realization of 
$ S_{q,\xi}^0 $ in terms of elements of $ \A.$ Then it turns out that 
this realization admits an infinitesimal characterization of 
$ S_{q,\xi}^0$. A representation of $ S_{q,\xi}^0 $ in terms of a 
$\U$-covariant oscillator is also given. 

  Analogous to the example of Podle\'s $q$-sphere, embedding of 
$ S_{q,\xi}^0 $ in $ \A $ may be done by realizing its generators in 
terms of entries of $ T^2(0)$ matrix. 
Denoting $ T^2(0) $ by $ T, $ it is straightforward to verify that 
the embedding is given by
\bea
  & & Y_2 = g_1\, T_{2,2} + g_2\, T_{0,2} + g_3\, T_{-2,2}, \nn \\
  & & Y_1 = g_1\, T_{2,1} + g_2\, T_{0,1} + g_3\, T_{-2,1}, \nn \\
  & & Y_0 = g_1\, T_{2,0} + g_2\, T_{0,0} + g_3\, T_{-2,0}, 
\label{embedding1} \\
  & & Y_{-1} = g_1\, T_{2,-1} + g_2\, T_{0,-1} + g_3\, T_{-2,-1}, \nn \\
  & & Y_{-2} = g_1\, T_{2,-2} + g_2\, T_{0,-2} + g_3\, T_{-2,-2}, \nn 
\eea
where the coefficients $ g_1, g_2 $ and $ g_3 $ need to satisfy the 
constraint
\beq
  g_1 g_3 = \frac{[3]!}{[4]} g_2^2.   \label{constraint1}
\eeq
In this embedding, the radius $r$ and the parameter $ \xi $ are given as 
functions of $g_2$:
\beq
   r = ([2]g_2)^2, \qquad 
   \xi = \frac{[6]}{[3]} g_2.    \label{rxi1}
\eeq
This embedding allows us to treat $ S_{q,\xi}^0 $ as a subalgebra of 
$ \A.$ This fact suggests that $ S_{q,\xi}^0 $ has an infinitesimal 
characterization {\it \`a la} Koornwinder \ocite{Koo}. To demonstrate 
this, we extend the left and right actions of $U_q[su(2)]$ on $SU_q(2)$ 
defined in Ref. \cite{Koo} to supersymmetric case. In the followings we 
assume $ u, v \in \U;\, a, b \in \A $, and use Sweedler's notation to 
denote coproducts: {\it e.g.} $\Delta(a) = \sum a_{(1)} \otimes a_{(2)}$. 
With a slight change of notation from Ref. \cite{Koo}, we now define 
elements $ u \odot a $ and $ a \odot u $ of $ \A $ by
\bea
  & & u \odot a = ({\rm id} \otimes u)(\Delta(a)) = \sum (-1)^{\h{u}\h{a}_{(1)}} a_{(1)} \pair{u}{a_{(2)}},
  \label{ua} \\
  & & a \odot u = (-1)^{\h{a}\h{u}}(u \otimes {\rm id})(\Delta(a)) = \sum (-1)^{\h{a}\h{u}}
      \pair{u}{a_{(1)}} a_{(2)}.
  \label{au}
\eea
The coassociativity of $\A$ then leads to
\beq
  (uv) \odot a = u \odot (v \odot a), \qquad 
  a \odot (uv) = (a \odot u) \odot v.
  \label{uva}
\eeq
Moreover, we also have
\bea
   & & u \odot (ab) = \sum (-1)^{\h{u}_{(2)}\h{a}} (u_{(1)} \odot a) (u_{(2)} \odot b),
   \nn \\ 
   & & (ab) \odot u = \sum (-1)^{\h{u}_{(1)}\h{b}} (a \odot u_{(1)})(b \odot u_{(2)}).
   \label{uab}
\eea
Thus, $ u \odot a $ and $ a \odot u$ define left and right actions of $u$ on $a$, respectively. 
The actions of generators of $ \U $ on $T^{\ell}$-matrix of $ \A $ are 
calculated by using (\ref{RepU}). Explicitly, the left actions are given by
\bea
   & & K^{\pm 1}\odot\T{\ell}{m_1m_2}(\lambda) = q^{\pm m_2/2}\, \T{\ell}{m_1m_2}(\lambda),
   \nn \\
   & & v_+ \odot \T{\ell}{m_1m_2}(\lambda) = (-1)^{\ell+m_1+\lambda}
       \sqrt{[\ell-m_2][\ell+m_2+1] \varrho} \,\T{\ell}{m_1\,m_2+1}(\lambda),
   \label{genT} \\
   & & v_- \odot \T{\ell}{m_1m_2}(\lambda) = (-1)^{m_1+m_2+\lambda+1} 
       \sqrt{[\ell+m_2][\ell-m_2+]\varrho}\, \T{\ell}{m_1\,m_2-1}(\lambda),
   \nn
\eea
while the right actions read
\bea
   & & \T{\ell}{m_1m_2}(\lambda) \odot K^{\pm 1} = q^{\pm m_1/2}\, \T{\ell}{m_1m_2}(\lambda),
   \nn \\
   & & \T{\ell}{m_1m_2}(\lambda) \odot v_+ = (-1)^{\ell+m_2+\lambda} 
       \sqrt{[\ell-m_1+1][\ell+m_1]\varrho}\, \T{\ell}{m_1-1\,m_2}(\lambda),
       \label{Tgen} \\
   & & \T{\ell}{m_1m_2}(\lambda) \odot v_- = (-1)^{m_1+m_2+\lambda} 
       \sqrt{[\ell+m_1+1][\ell-m_1]\varrho}\, \T{\ell}{m_1+1\,m_2}(\lambda).
   \nn
\eea

   An element $u \in \U $ possessing a coproduct structure 
$ \Delta(u) = g \otimes u + u \otimes g^{-1}$ with $ g \in \U $ being 
a group-like element, is known as \textit{twisted primitive} with respect 
to $g$. For a twisted primitive element $ u, $ it is straightforward to 
verify that 
\bea
   & & u \odot a = 0 \quad {\rm and} \quad u \odot b = 0 \quad \Rightarrow \quad u \odot (ab) = 0,
   \label{uaubuab} \\
   & & a \odot u = 0 \quad {\rm and} \quad b \odot u = 0 \quad \Rightarrow \quad (ab) \odot u = 0.
   \label{aubuabu}
\eea
Thus a set of elements of $ \A $ annihilated by a twisted primitive element $u$ form 
a subalgebra of $\A.$ Indeed, the quantum supersphere $ S^{0}_{q,\xi} $ embedded into $\A $ 
is a subalgebra of $ \A $ that is annihilated by the twisted primitive element 
$ {\cal P}_R$
\bea
   & & {\cal P}_R = -\sqrt{g_3}\, v_+ + \sqrt{g_1}\, v_-,
   \label{PR} \\
   & & Y_k \odot {\cal P}_R = 0, \qquad k = \pm 2,\; \pm 1,\; 0. \label{YkPR} 
\eea
The algebra $\U$ has three twisted primitive elements: 
$ K-K^{-1}, v_+$ and $ v_-. $ 
However, $ {\cal P}_R $ consists of only odd twisted primitive elements.  
This is a difference from the $q$-sphere for $SU_q(2)$. In that example, 
all the twisted primitive elements contribute to the annihilation 
operator of $q$-sphere. 

  We now turn to an oscillator realization of $ S^{0}_{q,\xi}. $ 
In Ref. \cite{TW}, a $\U$-covariant oscillator algebra is introduced. 
This oscillator algebra is generated by a pair of even creation/annihilation 
operators $ (\oab,\oa)$, and an odd operator $\oc$ obeying the relations 
\bea
  & & \oab \oc = q \oc \oab, \qquad
      \oa \oc = q^{-1} \oc \oa, \qquad
      \oa \oab -q^{-2} \oab \oa = 1,
  \nn \\
  & & \oc^2 = q^{-1}[2] \oab \oa + \varpi^{-1}. \label{osposc}
\eea
These relations are determined by two steps. First, the action of $\U$ on 
the oscillator is defined via the coproduct of $ \U.$ Then the commutation 
properties of the oscillator is fixed by demanding that the triplet 
$ (\oab,\oc,-\oa) $ transforms under the $ \ell = 1 $ representation of 
$\U.$ This suggests that the $\U$-covariant oscillator has a close kinship 
to the quantum space discussed in Sec.\ref{Acovariant}\ref{l10plane} 
Indeed, the $\U$-covariant oscillator is isomorphic to the quantum plane 
for a special value of $r$:
\beq
  z_1 = \oab, \qquad z_0 = \oc, \qquad z_{-1} = -\oa,
  \qquad r = -q \varpi^{-1}.
  \label{zosc}
\eeq
Employing the $\U$-covariant oscillator, it is possible to realize 
$ S^{0}_{q,\xi}$: 
\bea
  & & Y_2 = \oab^2, \qquad 
      Y_1 = q^{-1/2} \left( \frac{[4]}{[2]} \right)^{1/2} \oab \oc, 
  \nn \\
  & & Y_0 = -\frac{\sqrt{[4]!}}{q[2]} \oab \oa - \frac{q^{-1/2}}{\varpi} 
      \left( \frac{[4]}{[3]!} \right)^{1/2},
  \label{Yosc} \\
  & & Y_{-1} = -q^{-1/2} \left( \frac{[4]}{[2]} \right)^{1/2} \oc \oa, 
  \qquad 
  Y_{-2} = \oa^2. \nn 
\eea
In this realization, the radius $r$ and the parameter $\xi$ of 
$ S^{0}_{q,\xi} $ assume the following values
\beq
  r = \frac{q^2}{\varpi^2}\frac{[4]}{[3]!}, \qquad
  \xi = \frac{[6]}{[3]!} r^{1/2}.
  \label{rxiosc}
\eeq
An advantage of this realization is that we may represent 
$ S^{0}_{q,\xi} $ with matrices, 
since matrix representation of $\U$-covariant oscillator via that of 
Biedenharn-Macfarlane $q$-boson algebra exists\ocite{TW}.

%
\section{CONCLUDING REMARKS}
\label{CRemarks}

  We have developed the general prescription for constructing $\A$-covariant 
algebras. By the method, four $\A$-covariant algebras have been obtained, 
namely, three quantum superspaces and 
a one-parameter family of quantum superspheres. 
The special cases of the quantum superspaces correspond to the Manin's quantum 
superplane and $\U$-covariant oscillator algebra. 
The quantum superspheres are realized by $\A$ so that it can be regarded 
as a subalgebra of $\A$. This subalgebra is characterized by the fact that it  
is annihilated by the right action of a particular combination 
of the twisted primitive elements of $\U.$ These are the similarities to 
the $q$-spheres for $ SU_q(2).$  It has also been shown that the 
quantum superspheres have $\U$-covariant oscillator realization that 
allows us to have matrix representations of the quantum supersphere. 

   We believe that the results of this paper are useful for making 
progress in constructing supersymmetric versions of noncommutative 
geometry. For instance, we may consider differential calculi on the 
quantum supersphere, then compute its curvature, metric and so on 
based on the framework of Ref. \cite{AC}. In connection with the classification 
of differential calculi, it is interesting to determine the dual coalgebra of 
 $ S_{q,\xi}^0. $ The corresponding computation for $ SU_q(2)$ $q$-sphere 
was made recently and applied to the classification of differential 
calculi\ocite{HK}. 
Furthermore, the general method in Sec. \ref{Acovariant}A is applied 
to construct higher dimensional quantum superspaces by taking higher 
values of $\ell. $ 

  It is worth pointing out that, because of the similarity of the 
representation theory of $ \U $ to that of $ U_q[su(2)]$, 
the method developed in Sec. \ref{Acovariant}A is valid for 
$ SU_q(2). $ This procedure allows us to treat quantum plane, 
deformed oscillator and Podle\'s $q$-sphere in a unified way. 
The reason for $q$-spheres being a one-parameter family becomes 
clear in this framework. It is also possible to construct hitherto 
unknown higher dimensional $ SU_q(2)$-covariant quantum spaces. 
We will present these results elsewhere. 

%
\section*{ACKNOWLEDGEMENTS}

  The work of N.A. is partially supported by the grants-in-aid 
from JSPS, Japan (Contract No. 15540132). The work of the other author
(R.C.) is partially supported by the grant DAE/2001/37/12/BRNS, 
Government of India.  

%
\setcounter{section}{0}
\renewcommand{\thesection}{Appendix \Alph{section}}
\renewcommand{\thesubsection}{\Alph{section}.\arabic{subsection}}
\renewcommand{\theequation}{\Alph{section}.\arabic{equation}}
\setcounter{equation}{0}
\section{CGC OF $\U$}
\label{CGCderivation}

  In this appendix, the general expression of CGC in our conventions is 
derived. Let us write the highest weight vector as
\beq
  \e{\ell}{\ell}(\ell_1,\ell_2,\Lambda) = \sum_{m_1,m_2} A_{m_1m_2} \, 
  \e{\ell_1}{m_1}(\lambda) \otimes \e{\ell_2}{m_2}(\lambda),
  \label{ehighest}
\eeq
where, for simplicity of notation, the CGC is denoted by $ A_{m_1 m_2}$. 
The highest weight condition 
$ \Delta(v_+)\, \e{\ell}{\ell}(\ell_1,\ell_2,\Lambda) = 0 $ 
gives the recurrence relation for $ A_{m_1m_2}$
\bea
  & & \sqrt{[\ell_1-m_1][\ell_1+m_1+1]} q^{-m_2/2} A_{m_1m_2} \nn \\
  & & \hspace{2cm}
  - (-1)^{\ell_1-m_1+\lambda} q^{(m_1+1)/2} \sqrt{[\ell_2+m_2][\ell_2-m_2+1]}
  A_{m_1+1\, m_2-1} = 0.
  \label{rechigh}
\eea
It is easy to find the solution of this relation
\bea
  A_{m_1m_2} &=& (-1)^{\lambda(\ell_1-m_1) + (\ell_1-m_1)(\ell_1-m_1+1)/2}
  q^{(\ell+1)(\ell_1-m_1)/2}
  \nn \\
  &\times& \left(
    \frac{[\ell_1+\ell_2-\ell]!\,[\ell_1+m_1]!\,[\ell_2+m_2]!}
         {[\ell_2-\ell_1+\ell]!\,[2\ell_1]!\,[\ell_1-m_1]!\,[\ell_2-m_2]!}
  \right)^{1/2} A_{\ell_1\, \ell-\ell_1}.
  \label{Am1m2}
\eea
The normalization of the highest weight vector (\ref{ehighest})
determines $ A_{\ell_1\,\ell-\ell_1} $ as follows: 
\bea
  & & (\e{\ell}{\ell}(\ell_1,\ell_2,\Lambda),\e{\ell}{\ell}(\ell_1,\ell_2,\Lambda))
   = \sum_{m_1,m_2} (-1)^{(\ell_1-m_1+\lambda)(\ell_2-m_2+\lambda)} A_{m_1m_2}^2
  \nn \\
  & & = (-1)^{(\ell_1+\lambda)(\ell_2+\ell+\lambda)} 
      \frac{[\ell_1+\ell_2-\ell]!}{[\ell_2-\ell_1+\ell]!\, [2\ell_1]!} A_{\ell_1\, \ell-\ell_1}^2
  \nn \\
  & & \times \sum_{m_1} (-1)^{m_1(\ell_1+\ell_2+\ell+1)} q^{(\ell+1)(\ell_1-m_1)} 
      \frac{[\ell_1+m_1]!\, [\ell_2+\ell-m_1]!}{[\ell_1 -m_1]!\, [\ell_2-\ell+m]!}.
      \nn
\eea
The summation over $m_1$ is computed by using the formula (\ref{sumbnc}). 
Setting
\[
   a = \ell_1-m_1, \quad k= \ell_1+\ell_2-\ell, \quad
   -n = \ell_1-\ell_2+\ell+1, \quad -r = -\ell_1 + \ell_2 + \ell+1
\]
in (\ref{sumbnc}), while noticing that all these quantities are positive 
integers, we obtain
\bea
  & & \sum_{m_1} (-1)^{(\ell_1-m_1)(\ell_1+\ell_2+\ell+1)} q^{(\ell+1)(\ell_1-m_1)} 
      \frac{[\ell_1+m_1]!\, [\ell_2+\ell-m_1]!}{[\ell_1 -m_1]!\, [\ell_2-\ell+m]!}
  \nn \\
  & & \qquad = q^{(-\ell_1+\ell_2+\ell+1)(\ell_1+\ell_2-1)/2}
      \frac{[\ell_1+\ell_2+\ell+1]!\, [\ell_1-\ell_2+\ell]! \, [-\ell_1+\ell_2+\ell]!}
           {[\ell_1+\ell_2-\ell]!\, [2\ell+1]!}.
  \nn 
\eea
Thus the norm of the highest weight vector reads
\bea
 & & ||\e{\ell}{\ell}(\ell_1,\ell_2,\Lambda)||^2 \nn \\
 & & \qquad = (-1)^{\lambda(\ell_1+\ell_2+\ell+1)}
 q^{(-\ell_1+\ell_2+\ell+1)(\ell_1+\ell_2-1)/2} A_{\ell_1\, \ell-\ell_1}^2 
 \frac{[\ell_1+\ell_2+\ell+1]!\,[\ell_1-\ell_2+\ell]!}{[2\ell+1]!\,
 [2\ell_1]!}.\nn
\eea
This leads to
\beq
  A_{\ell_1\, \ell-\ell_1} =  q^{-(-\ell_1+\ell_2+\ell+1)(\ell_1+\ell_2-1)/4}
  \left(\frac{[2\ell+1]!\,[2\ell_1]!}{[\ell_1+\ell_2+\ell+1]!\,
  [\ell_1-\ell_2+\ell]!}\right)^{1/2},
  \label{Adet}
\eeq
where the phase is chosen such that the expression coincides the result in 
Ref. \cite{MM}. 

  To obtain the other vectors in this irreducible representation  
we need the following results, which may be verified by induction:
\bea
  & & v_-^k \e{\ell}{m}(\lambda) = (-1)^{(\ell-m)k + k(k+1)/2} 
      \left(
         \frac{[\ell+m]!\,[\ell-m+k]!}{[\ell-m]!\,[\ell+m-k]!} \varrho^k
      \right)^{1/2} \e{\ell}{m-k}(\lambda),
   \label{vmk} \\
  & & \Delta(v_-^n) = \sum_{k=0}^n \bn{n}{k} (-1)^{k(n-k)} v_-^{n-k}K^k \otimes v_-^k K^{-n+k}.
   \label{Dvmn}
\eea
Using (\ref{vmk}) we obtain
\beq
  \e{\ell}{m}(\ell_1,\ell_2,\Lambda) = (-1)^{(\ell-m)(\ell-m+1)/2}
  \left(
     \frac{[\ell+m]!}{[2\ell]!\,[\ell-m]! \,\varrho^{\ell-m}}
  \right)^{1/2}
  \Delta(v_-^{\ell-m})\, \e{\ell}{\ell}(\ell_1,\ell_2,\Lambda).
\eeq
The right hand side is computed by using (\ref{Dvmn}) and (\ref{ehighest}). 
After some algebra, we derive
\bea
  & & \Delta(v_-^{\ell-m})\, \e{\ell_1}{m_1}(\lambda) \otimes \e{\ell_2}{m_2}(\lambda)
  \nn \\
  & & = \sum_k \bn{\ell-m}{k} (-1)^{k(\ell_2-m_2+\lambda) + (\ell_1-m_1)(\ell-m)
        + (\ell-m)(\ell-m+1)/2}
      q^{km_1/2 + (-\ell+m+k)m_2/2}
  \nn \\
  & & \quad \times
  \left(
    \frac{[\ell_1+m_1]!\,[\ell_2+m_2]!\,[\ell_1+\ell-m_1-m-k]!\,[\ell_2-m_2+k]!}
    {[\ell_1-m_1]!\,[\ell_2-m_2]!\,[\ell_1-\ell+m_1+m+k]!\,[\ell_2+m_2-k]!}
     \varrho^{\ell-m}
  \right)^{1/2}
   \nn \\
  & & \quad \times \;
      \e{\ell_1}{m_1-\ell+m+k}(\lambda) \otimes \e{\ell_2}{m_2-k}(\lambda).
  \label{vmlmonvec}
\eea
Equations (\ref{Am1m2}) and (\ref{Adet}) allow us to derive 
\bea
  & & \e{\ell}{m}(\ell_1,\ell_2,\Lambda)  
  \nn \\
  & & \qquad = \sum_{m_1,m_2} 
     (-1)^{(\ell_1-\ell+m_2)(\ell-m+\lambda) + (\ell_1-\ell+m2)(\ell_1-\ell+m_2+1)/2}
  \nn \\
  & & \qquad \times q^{m_2(m+1)/2 + (\ell_1-\ell_2)(\ell_1+\ell_2+1)/4-\ell(\ell+1)/4}
  \nn \\
  & & \qquad \times \left(
    [2\ell+1] \frac{[\ell_1+\ell_2-\ell]!\,[\ell+m]!\,[\ell-m]!\,[\ell_1-m_1]!\,[\ell_2-m_2]!}
                  {[\ell_1+\ell_2+\ell+1]!\,[\ell_1-\ell_2+\ell]!\,[-\ell_1+\ell_2+\ell]!\,
                  [\ell_1+m_1]!\,[\ell_2+m_2]!}
   \right)^{1/2}
  \nn\\
  & & \qquad \times
    \sum_{k}(-1)^{k(k-1)/2+ k(\ell_1+\ell_2-m)} q^{k(\ell+m+1)/2}
  \nn \\
  & & \qquad \times
      \frac{[\ell_1+\ell-m_2-k]!\,[\ell_2+m_2+k]!}
      {[k]!\,[\ell-m-k]!\,[\ell_1-\ell+m_2+k]!\,[\ell_2-m_2-k]!}
      \e{\ell_1}{m_1}(\lambda) \otimes \e{\ell_2}{m_2}(\lambda).
      \label{irrepdecomp}
\eea

%
\setcounter{equation}{0}
\section{SUMMATION FORMULA FOR KULISH SYMBOL}

  A summation formula for Kulish symbol which is used in the previous 
section is derived in this appendix. Corresponding binomial coefficient  
is usually defined for $ y \geq x \geq 0 $ as
\beq
  \bn{y}{x} = \frac{[y]!}{[y-x]![x]!} = \frac{[y][y-1] \cdots [y-x+1]}{[x]!}.
  \label{bnc}
\eeq
We may extend this to $ y < 0 $ by using the property 
\beq
  [-n] = (-1)^{n+1}[n], \qquad n \geq 0, \label{negativeK}
\eeq
and define the binomial coefficient for negative $y$ by the right most 
formula of (\ref{bnc}):
\bea
  \bn{y}{x} &=& \frac{[-(-y)][-(-y+1)] \cdots [-(-y+x-1)]}{[x]!} \nn \\
            &=& (-1)^{xy+x(x+1)/2} \frac{[-y][-y+1] \cdots [-y+x-1]}{[x]!}. \nn 
\eea
Thus for $ y < 0 $ and $ x \geq 0, $ the binomial coefficient is given by
\beq
  \bn{y}{x} =  (-1)^{xy+x(x+1)/2} \bn{x-y-1}{x}. \label{negativebnc}
\eeq
The following formula, proved by induction, is found in Ref. \cite{MM}:
\beq
  \bn{n+r}{k} = \sum_{a} \bn{n}{k-a} \bn{r}{a} (-1)^{(k-a)(r-a)}
   q^{(r-a)(n+r)/2-r(n-k+r)/2},
   \label{sumbn}
\eeq
where $ a $ runs over any positive integers such that the arguments of $ [x] $ in 
the binomial coefficients are non-negative. 
We assume $ n, r < 0 $, and apply (\ref{negativebnc}) to this formula
\bea
  & & {\rm l.h.s.} = (-1)^{k(n+r)+k(k+1)/2} \bn{k-n-r-1}{k}, \nn \\
  & & {\rm r.h.s.} = (-1)^{k(n+r)+k(k+1)/2} q^{rk/2} \sum_a (-1)^{an} 
      q^{-a(n+r)/2} \bn{k-a-n-1}{k-a} \bn{a-r-1}{a}.\nn
\eea
From this, the following summation formula of Kulish symbols is obtained
\bea
  & & \sum_a (-1)^{an} q^{-a(n+r)/2} 
  \frac{[k-a-n-1]!\, [a-r-1]!}{[k-a]!\, [a]!} 
  \nn \\
  & & \hspace{3cm} = 
  q^{-rk/2} \frac{[k-n-r-1]!\, [-n-1]!\, [-r-1]!}{[k]!\, [-n-r-1]!}.
  \label{sumbnc}
\eea

\newpage
%
\setcounter{equation}{0}
\section{CGC FOR $ \bf 1 \otimes 1 = 2 \oplus 1 \oplus 0 $}
\label{L1CGC}

  The following tables contain the values of 
$ C^{\;1 \quad \ 1 \quad \ \;\ell}_{m_1\ m-m_1\ m} $ 
for a given $ \ell. $ The columns provide the values of 
$m_1$, while the rows indicate $ m.$ 
The rightmost column of Table~\ref{C12}, titled as "OF", indicates the overall factors 
that are common to all entries in the row. In Table~\ref{C11} and Table~\ref{C10}, 
the overall factors (OF) are common for all entries of the tables.

\begin{table}[ht]
{\small
 \begin{center}
 \begin{tabular}{c||c|c|c||c}
    & 1 & 0 & $-1$ & OF \\ \hline
  2 & 1 & 0 &  0 & 1\\ \hline
  1 & $(-1)^{\lambda} q^{1/2}$ & $q^{-1/2}$ & 0  & 
    $ {\displaystyle \left( \frac{[2]}{[4]}\right)^{1/2}} $ \\ \hline
  0 & $q$ & $(-1)^{\lambda} [2]$ & $q^{-1}$ &
    $ {\displaystyle \frac{[2]}{\sqrt{[4]!}} }$ \\ \hline
  $-1$ & 0 & $q^{1/2}$ & $(-1)^{\lambda} q^{-1/2}$ & 
    ${\displaystyle \left( \frac{[2]}{[4]}\right)^{1/2}}$ \\ \hline
  $-2$ & 0 & 0 & 1 & 1
 \end{tabular}
 \caption{$ \ell=2$} \label{C12}
 \end{center}
 }
\end{table}


\begin{table}[ht]
{\small
 \begin{center}
 \begin{tabular}{cc}
  \begin{minipage}{0.5\hsize}
  \begin{center}
   \begin{tabular}{c||c|c|c}
      & 1 & 0 & $-1$  \\ \hline
    1 & $q^{-1/2}$ & $(-1)^{\lambda+1} q^{1/2}$ & 0  \\ \hline
    0 & $(-1)^{\lambda+1}$ & $q^{1/2}+q^{-1/2}$ & $(-1)^{\lambda}$ \\ \hline
    $-1$ & 0 & $(-1)^{\lambda+1} q^{-1/2}$ & $q^{1/2}$ 
   \end{tabular}
 \caption{$\ell = 1$,  $ {\rm OF} = {\displaystyle \left( \frac{[2]}{[4]}\right)^{1/2}}.$
 \label{C11}
}
 \end{center}
 \end{minipage}
 \begin{minipage}{0.5\hsize}
 \begin{center}
   \begin{tabular}{c||c|c|c}
       & 1 & 0 & $-1$ \\ \hline
     0 & $q^{-1/2}$ & $(-1)^{\lambda+1}$ & $-q^{1/2}$ 
   \end{tabular}
  \caption{$\ell=0$,  ${\rm OF} = {\displaystyle \frac{1}{\sqrt{[3]}} }$.}
  \label{C10}
  \end{center}
  \end{minipage}
 \end{tabular}
 \end{center}
}
\end{table}



%
\section{CGC FOR $ \bf 2 \otimes 2 = 4 \oplus 3 \oplus 2 \oplus 1 \oplus 0 $}
\label{L2CGC}

  The following tables contain the values of 
$ C^{\;2 \quad \ 2 \quad \ \;\ell}_{m_1\ m-m_1\ m} $ 
for a given $ \ell. $ The columns provide the values of 
$m_1$, while the rows indicate $ m.$ 
The right most column, titled as "OF", indicates the overall factors 
common to all entries in the row. 


\begin{landscape}

\begin{table}[ht]
{\small
\begin{center}
\begin{tabular}{c||c|c|c|c|c||c}
    &  2 & 1 &  0 & $-1$ & $-2$ & {\rm OF}\\
  \hline \hline
  4 & 1 & 0 & 0 & 0 & 0 & 1\\
  \hline
  3 &
  $\displaystyle (-1)^{\lambda} q$  &
  $\displaystyle  q^{-1}$ &
  0 & 0 & 0 & 
  $\displaystyle \left( \frac{[4]}{[8]} \right)^{1/2}$
  \\
  \hline
  2 & 
  $\displaystyle q^2$ & 
  $\displaystyle (-1)^{\lambda} \left( \frac{[2][4]}{[3]} \right)^{1/2}$ &
  $\displaystyle q^{-2}$ & 
  0 &
  0 & 
  $\displaystyle \left( \frac{[3][4]}{[7][8]} \right)^{1/2}$ 
  \\
  \hline
  1 &
  $\displaystyle (-1)^{\lambda} q^3 \sqrt{[4]!}$  & 
  $\displaystyle q [4][3]$  &
  $\displaystyle (-1)^{\lambda} q^{-1} [4][3]$  &
  $\displaystyle q^{-3} \sqrt{[4]!}$ &
  0 &
  $\displaystyle \left( \frac{[5]!}{[8]!} \right)^{1/2}$
  \\
  \hline
  0 &
  $\displaystyle q^4$ &
  $\displaystyle (-1)^{\lambda}q^2 [4]$ &
  $\displaystyle [4][3]$ &
  $\displaystyle (-1)^{\lambda} q^{-2} [4]$  &
  $\displaystyle q^{-4}$ &
  $\displaystyle \frac{[4]!}{\sqrt{[8]!}}$
  \\
  \hline
  $-1$ &
  0 &
  $\displaystyle q^3 \sqrt{[4]!}$  &
  $\displaystyle (-1)^{\lambda} q [4][3]$  &
  $\displaystyle q^{-1} [4][3]$ &
  $\displaystyle (-1)^{\lambda} q^{-3} \sqrt{[4]!}$ & 
  $\displaystyle \left( \frac{[5]!}{[8]!} \right)^{1/2}$
  \\
  \hline
  $-2$ &
  0 & 0 & 
  $\displaystyle q^2$ & 
  $\displaystyle (-1)^{\lambda} \left(\frac{[2][4]}{[3]}\right)^{1/2}$ & 
  $\displaystyle q^{-2}$ &
  $\displaystyle \left( \frac{[3][4]}{[7][8]} \right)^{1/2}$ 
  \\
  \hline
  $-3$ &
  0 & 0 & 0 &
  $\displaystyle q $ &
  $\displaystyle (-1)^{\lambda} q^{-1}$  & 
  $\displaystyle \left( \frac{[4]}{[8]} \right)^{1/2}$
  \\
  \hline
  $-4$ & 0 & 0 & 0 & 0 & 1 & 1
\end{tabular}
\caption{$ \ell = 4$}
\end{center}
}
\end{table}

%
\begin{table}[ht]
{\small
\begin{center}
\begin{tabular}{c||c|c|c|c|c||c}
   & $ 2$ & $1$ & $0$ & $-1$ & $-2$ & {\rm OF} \\
  \hline \hline
   3 & 
   $\displaystyle q^{-1}$ & 
   $\displaystyle (-1)^{\lambda+1} q$  & 
   0 & 
   0 & 
   0 & 
   $\displaystyle \left(\frac{[4]}{[8]}\right)^{1/2} $
  \\
  \hline
  2 &
  $\displaystyle (-1)^{\lambda+1}$  &
  $\mu_{31}$ &
  $\displaystyle (-1)^{\lambda}$  &
  0 &
  0 & 
  $\displaystyle \left( \frac{[4]!}{[6][8]} \right)^{1/2}$
  \\
  \hline
  1 &
  $q$ &
  $(-1)^{\lambda+1} \mu_{32}$ &
  $\mu_{33}$ &
  $(-1)^{\lambda+1} q^{-1}$ &
  0 &
  $\displaystyle \left( \frac{[7][3]}{[8]!} \right)^{1/2} [4]!$
  \\
  \hline
  0 &
  $(-1)^{\lambda+1} q^2$ &
  $q^{1/2} ([3] + q^2)$ &
  $(-1)^{\lambda} [3]\omega$ &
  $q^{-1/2} ([3] + q^{-2})$ &
  $(-1)^{\lambda} q^{-2}$ &
  $\displaystyle \left( \frac{[7]}{[8]!} \right)^{1/2}[4]!$
  \\
  \hline
  $-1$ &
  0 &
  $(-1)^{\lambda+1} q$ &
  $\mu_{32}$ &
  $(-1)^{\lambda+1} \mu_{33}$ &
  $q^{-1}$ &
  $\displaystyle \left( \frac{[7][3]}{[8]!} \right)^{1/2} [4]!$
  \\
  \hline
  $-2$ &
  0  &
  0 &
  $(-1)^{\lambda+1}$  &
  $\mu_{31}$ &
  $(-1)^{\lambda}$ & 
  $\displaystyle \left( \frac{[4]!}{[6][8]} \right)^{1/2}$
  \\
  \hline
   $-3$ & 
   0 & 
   0 & 
   0 & 
   $(-1)^{\lambda+1} q^{-1}$ & 
   $q$ & 
   $\displaystyle \left(\frac{[4]}{[8]}\right)^{1/2}$ 
\end{tabular}

\[
 \mu_{31} = \frac{\sqrt{ [4]! }}{[2]^2}\omega,\qquad
 \mu_{32} = \left( \frac{[4]}{[3]!} \right)^{1/2} (q^2+q^{-1/2}[2]), \qquad
 \mu_{33} = \left( \frac{[4]}{[3]!} \right)^{1/2} (q^{-2}-q^{1/2}[2]).
\]
\caption{$\ell = 3$}
\end{center}
}
\end{table}


\begin{table}[ht]
{\small
\begin{center}
\begin{tabular}{c||c|c|c|c|c||c}
   & $2$ & $ 1$ & $ 0$ & $-1$ & $-2$ & {\rm OF} \\
  \hline \hline
   2 & 
   $q^{-3/2}$  & 
   $\displaystyle (-1)^{\lambda+1} \left( \frac{[3]!}{[4]} \right)^{1/2}$ & 
   $-q^{3/2}$  & 
   0  & 
   0 &
   $\displaystyle  \left( \frac{[3][4]}{[6][7]} \right)^{1/2}$
  \\
  \hline
  1 &
  $\displaystyle (-1)^{\lambda} q^{-1/2} \left( \frac{[3]!}{[4]} \right)^{1/2}$  &
  $q^{-1/2} \mu_{21}$  &
  $(-1)^{\lambda+1} q^{1/2} \mu_{22}$ &
  $\displaystyle -q^{1/2} \left( \frac{[3]!}{[4]} \right)^{1/2}$ &
  0  &
  $\displaystyle \left( \frac{[3][4]}{[6][7]} \right)^{1/2}$
  \\
  \hline
  0 & 
  $q^{1/2}$ &
  $(-1)^{\lambda+1} \mu_{22}$ &
  $\mu_{23}$ &
  $(-1)^{\lambda+1} \mu_{21}$ &
  $-q^{-1/2}$ &
  $\displaystyle \left( \frac{[3][4]}{[6][7]} \right)^{1/2}$
  \\
  \hline
  $-1$ &
  0 &
  $\displaystyle q^{-1/2} \left( \frac{[3]!}{[4]} \right)^{1/2}$ &
  $(-1)^{\lambda} q^{-1/2} \mu_{21}$ &
  $-q^{1/2} \mu_{22}$ & 
  $\displaystyle (-1)^{\lambda+1} q^{1/2} \left( \frac{[3]!}{[4]} \right)^{1/2}$ &
  $\displaystyle \left( \frac{[3][4]}{[6][7]} \right)^{1/2}$
  \\
  \hline
   $-2$ & 
   0 &
   0 &
   $q^{-3/2}$  & 
   $\displaystyle (-1)^{\lambda+1} \left( \frac{[3]!}{[4]} \right)^{1/2}$ & 
   $-q^{3/2}$  & 
   $\displaystyle  \left( \frac{[3][4]}{[6][7]} \right)^{1/2}$
\end{tabular}

\[
  \mu_{21} = \frac{[2]}{[4]} (q^{-2} + q^{3/2}[2]), \qquad
  \mu_{22} = \frac{[2]}{[4]} (q^2-q^{-3/2}[2]), \qquad
  \mu_{23} = \frac{[2]^2}{[4]} \left( \frac{[8]}{[4]} + 2\frac{[4]}{[2]} + 1 \right).
\]
\end{center}
\caption{$ \ell = 2$}
}
\end{table}


\begin{table}[ht]
{\small
\begin{center}
\begin{tabular}{c||c|c|c|c|c||c}
    & $ 2$ & $1$ & $0$ & $-1$ & $-2$ & {\rm OF} \\
  \hline \hline
   1 & 
   $\displaystyle q^{-3/2}$ & 
   $\displaystyle (-1)^{\lambda+1} q^{-1/2}  \left( \frac{[3]!}{[4]} \right)^{1/2}$ & 
   $\displaystyle -q^{1/2} \left( \frac{[3]!}{[4]} \right)^{1/2}$ & 
   $\displaystyle (-1)^{\lambda} q^{3/2}$ & 
   0 &
   $\displaystyle \left( \frac{[3]!}{[5][6]} \right)^{1/2}$ 
  \\
  \hline
  0 &
  $(-1)^{\lambda+1} q^{-1/2}$  &
  $q^{-1} \mu_{11}$  &
  $\displaystyle (-1)^{\lambda+1}\frac{[3]}{[4]}\omega$  &
  $-q \mu_{12}$  &
  $(-1)^{\lambda+1}q^{1/2}$ &
  $\displaystyle \left( \frac{[3][4]}{[5][6]} \right)^{1/2} $
  \\
  \hline
  $-1$ &
  0 &
  $\displaystyle (-1)^{\lambda+1} q^{-3/2}$ &
  $\displaystyle q^{-1/2} \left( \frac{[3]!}{[4]} \right)^{1/2}$  &
  $\displaystyle (-1)^{\lambda} q^{1/2} \left( \frac{[3]!}{[4]} \right)^{1/2}$ &
  $\displaystyle -q^{3/2}$ &
  $\displaystyle  \left( \frac{[3]!}{[5][6]} \right)^{1/2}$
\end{tabular}

\[
  \mu_{11} = q + q^{-1} \frac{[2]}{[4]},\qquad
  \mu_{12} = q^{-1} + q \frac{[2]}{[4]}.
\]
\end{center}
\caption{$ \ell = 1$}
}
\end{table}


\begin{table}[ht]
{\small
\begin{center}
\begin{tabular}{c||c|c|c|c|c||c}
     & $2$ & $1$ & $0$ & $-1$ & $-2$ & {\rm OF} \\
  \hline \hline
   0 & 
   $\displaystyle q^{-1}$ & 
   $\displaystyle (-1)^{\lambda+1} q^{-1/2}$ & 
   $\displaystyle -1$ & 
   $\displaystyle (-1)^{\lambda} q^{1/2}$ & 
   $\displaystyle\ q $ &
   $\displaystyle \frac{1}{\sqrt{[5]}}$
\end{tabular}
\end{center}
\caption{$ \ell = 0$}
}
\end{table}

\end{landscape}

%
%

\end{document}